\documentclass[11pt]{article}
\usepackage{amsmath,epsfig}
\usepackage{amsthm}
\usepackage{amssymb}
\numberwithin{equation}{section}

\def\nm{\noalign{\medskip}}
\usepackage{geometry} 
\newcommand{\ds}{\displaystyle}
\usepackage{url}
\usepackage{hyperref}
\newtheorem{lemma}{Lemma}
\newtheorem{theorem}{Theorem}
\newtheorem{prop}{Proposition}
\newtheorem{asump}{Assumption}[section]
\newtheorem{rmk}{Remark}[section]

\newcommand{\R}{\mathbb{R}}
\def\bea{\begin{eqnarray}}  \def\eea{\end{eqnarray}}
 \def\p{\partial}
\def \Vh0{\stackrel{\circ}{V}_h} 
   
\def\Om{\Omega}  
   
  \def\f{\frac}  

\def\beas{\begin{eqnarray*}} \def\eeas{\end{eqnarray*}}

\title{Double-negative electromagnetic metamaterials due to chirality}

\author{
Habib Ammari\thanks{\footnotesize Department of Mathematics, 
ETH Z\"urich, 
R\"amistrasse 101, CH-8092 Z\"urich, Switzerland (habib.ammari@math.ethz.ch, wei.wu@sam.math.ethz.ch, sanghyeon.yu@sam.math.ethz.ch).} \and Wei Wu\footnotemark[2] \and Sanghyeon Yu\footnotemark[2]}

\date{}

\begin{document}
\maketitle

\begin{abstract}
The aim of this paper is to provide a mathematical theory for understanding the mechanism behind the double-negative refractive index phenomenon in chiral materials. The design of double-negative metamaterials generally requires the use of two different kinds of subwavelength resonators, which may limit the applicability of double-negative metamaterials. Herein, we rely on media that consist of only a single type of dielectric resonant element, and show how the chirality of the background medium induces double-negative refractive index metamaterial, which refracts waves negatively, hence acting as a superlens. Using plasmonic dielectric particles, it is proved that both the effective electric permittivity and the magnetic permeability can be negative near some resonant frequencies. A justification of the approximation of a plasmonic particle in a chiral medium by the sum of a resonant electric dipole and a resonant magnetic dipole, is provided. Moreover, the set of resonant frequencies is characterized. For an appropriate volume fraction of plasmonic particles with certain conditions on their configuration, a double-negative effective medium can be obtained when the frequency is near one of the resonant frequencies.

\end{abstract}

\def\keywords2{\vspace{.5em}{\textbf{  Mathematics Subject Classification
(MSC2000).}~\,\relax}}
\def\endkeywords2{\par}
\keywords2{35R30, 35C20.}

\def\keywords{\vspace{.5em}{\textbf{ Keywords.}~\,\relax}}
\def\endkeywords{\par}
\keywords{Plasmonic nanoparticles, sub-wavelength resonance, electric and magnetic resonant dipoles, chiral materials, effective medium theory, double-negative metamaterials.}

\section{Introduction}

The resolving power of conventional imaging systems is generally limited by the operating wavelength, which prevents imaging of subwavelength structures. However, systems having a negative refractive index can produce sharp images, with a potential to resolve subwavelength features \cite{pendry2000}.  
Negative index materials were first considered in \cite{VGV1}. Their recently added potential for subwavelength imaging led to an enormous interest in their properties \cite{smith, shalaev, soukoulis, veselago}. 

The negative sign of the refractive index arises in the description of electromagnetic properties of materials with simultaneously negative values of dielectric permittivity and magnetic permeability.  
A negative refractive index means that the phase velocity of a propagating wave  is opposite to the movement of the energy flux of the wave, represented by the Poynting vector \cite{veselago}.  

Negative refractive index materials are often referred to as double-negative metamaterials. 
The design of double-negative metamaterials generally requires the use of two different kinds of building blocks or specific subwavelength resonators presenting overlapping resonances. The negative index is caused by simultaneous resonant electric and magnetic responses.
The resonant electric dipole modes in the background medium contribute to the effective dielectric permittivity while the resonant magnetic dipole modes contribute to the effective magnetic permeability. Such a requirement limits the applicability of double-negative metamaterials.

Recently, it has been predicted and experimentally demonstrated  that introducing a single type of sub-wavelength resonant dielectric element into chiral materials leads to double-negative metamaterials \cite{JBP1, chiraleneg1,chiraleneg2}.

A material is defined to be chiral if it lacks any planes of mirror symmetry. In terms of electromagnetic responses, chiral materials exhibit a different refractive index for each polarization of the electromagnetic wave and are characterized by a cross coupling between the electric and magnetic dipoles along the same direction. The coupling strength is given by the magnitude of a quantity known as the chirality admittance which determines the bulk electromagnetic properties of chiral materials. Wave propagation inside chiral materials is investigated, for instance, in \cite{hamdache, laouadi,HA1,HA2,ima,m2as,martin,mitrea}. 

Plasmonic nanoparticles are sub-wavelength resonant dielectric elements. They exhibit quasi-static resonances, called plasmonic resonances. At or near these resonant frequencies, the plasmonic particles behave as strong electric dipoles. The plasmonic resonances are related to the spectra of the non-self adjoint Neumann-Poincar\'e type operators associated with the particle shapes. We refer the reader to \cite{pierre, matias1, HA4, sensing, hyeonbae1,hyeonbae2} for recent mathematical analysis of fundamental plasmonic resonance phenomena and their implications in subwavelength imaging.  

In this paper, we aim to understand the mechanism behind the double-negative refractive index phenomenon in chiral media. Herein, we rely on media that consists of plasmonic resonant particles, and show how to turn a chiral material into a negative refractive index metamaterial.
For this purpose, we first derive the leading-order perturbations in the electromagnetic fields in the far-field, which are caused by the presence of a plasmonic nanoparticle. To our knowledge, these asymptotic expansions, which are uniformly valid with respect to the frequency, have never been established. They generalize those derived in \cite{pierre, khelifi, vogelius} to chiral media. They show that the plasmonic nanoparticle can be approximated by the sum of resonant electric 
and magnetic dipole sources. Then, we completely characterize the set of resonant frequencies in terms of the plasmonic resonances of the nanoparticle in free space and the chiral admittance of the background medium. Finally, by using the point interaction approximation, 
we show that double-negative electromagnetic materials can be obtained by embedding in a chiral medium a large number of regularly spaced, randomly oriented plasmonic particles, each modeled as the sum of resonant electric and magnetic dipole sources. Near or at the resonant frequencies, the effective electric and magnetic properties of the medium can both be negative.   
 We recall that the idea of point interaction approximation goes back to Foldy's paper \cite{Foldy}. It is a natural tool to analyze a variety of interesting problems in the continuum limit. It was first applied  to the analysis of boundary value problems in regions with many small holes\cite{ozawa, ozawa2},  then in \cite{figari} to the heat conduction in material with many small holes, and in \cite{caflish} on sound propagation in bubbly fluid.

Our methodology in this paper follows the one introduced recently for achieving double-negative acoustic media using bubbles \cite{acoustic}. In acoustics, it is known that the air bubbles are subwavelength resonators \cite{Minnaert1933}. Due to the high contrast between the air density inside and outside an air bubble in a fluid,  a quasi-static acoustic resonance known as the Minnaert resonance occurs \cite{H3a}. At or near this resonant frequency, the size of a bubble can be up to three orders of magnitude smaller than the wavelength of the incident wave, and the bubble behaves as a strong monopole scatterer of sound. The Minnaert resonance phenomenon makes air bubbles good candidates for acoustic subwavelength resonators.
In \cite{acoustic}, it is proved that, using bubble dimers, the effective mass density and bulk modulus of the bubbly fluid can both be negative over a non empty range of frequencies. A bubble dimer consists of two identical separated bubbles.   It features two slightly different subwavelength resonances, called the hybridized Minnaert resonances. The hybridized Minnaert resonances are fundamentally different modes. One mode is a monopole as in the case of a single bubble, while the other one is a dipole. 
The resonance associated with the dipole mode is usually referred to as the anti-resonance. For an appropriate volume fraction, when the excitation frequency is close to the anti-resonance, a double-negative effective mass density and bulk modulus for bubbly media consisting of a large number of bubble dimers with certain conditions on their distribution are obtained. The dipole modes in the background medium contribute to the effective mass density while the monopole modes contribute to the effective bulk modulus.

The paper is organized as follows. In Section \ref{sec1}, we introduce some preliminaries on electromagnetic wave propagation in chiral materials. In Section \ref{sec2}, we derive an asymptotic expansion of the scattered electric and magnetic fields  by a plasmonic nanoparticles in a chiral material. We prove that the plasmonic nanoparticle   can be approximated as a pair of electric and magnetic dipole sources. We also characterize the set of resonant frequencies. 
In Section \ref{sec_EMT}, we derive a double-negative effective medium theory for plasmonic particles in chiral media near the resonant frequencies. In Appendix \ref{app_ball}, we provide some explicit calculations for the case of spherical chiral inclusions. 
  In Appendix \ref{subsec:layer_chiral}, we review layer potential formulations for electromagnetic scattering  in a chiral medium.

\section{Electromagnetic scattering in a chiral material} \label{sec1}

In this section, we consider an infinite chiral material in $\mathbb{R}^3$ with only one plasmonic particle. Let the particle $\Omega\subset\mathbb{R}^3$ be an open bounded set, with smooth boundary $\partial\Omega$. We assume that the particle $\Omega$ is centered at the origin and is of size $O(\delta)$, where $\delta >0$ is small. Denote $\Omega= \delta B$, $|B|=O(1)$. Obviously, $|\Omega|=\delta^3|B|=O(\delta^3)$, where $|\;\cdot \; |$ denotes the volume. The electric permittivity $\epsilon(x)$, permeability $\mu(x)$, and chiral admittance $\beta(x)$ in $\mathbb{R}^3$ satisfy
\begin{equation}\label{eq:parametersetting}
\beta = \left\{\begin{aligned}
&\beta_m\quad\mathrm{in}~\mathbb{R}^3 \setminus \overline{\Omega} \\
&0~~~\quad\mathrm{in}~\Omega
\end{aligned}
\right.,\quad
\epsilon = \left\{\begin{aligned}
&\epsilon_m~~\quad\mathrm{in}~\mathbb{R}^3 \setminus \overline{\Omega} \\
&\epsilon_c~~~\quad\mathrm{in}~\Omega
\end{aligned}\right.,\quad \mbox{ and } \mu(x) = \mu_m \mbox{  for  all } x \in \mathbb{R}^3. 
\end{equation}
Here, $\epsilon_m, \mu_m,$ and $\beta_m$ are positive constants and $\epsilon_c$ depends on the operating frequency $\omega$ and is assumed to be negative. 

\subsection{Drude-Born-Fedorov equations}

The starting points of this paper are the Maxwell equations and the constitutive relations for the chiral medium. Different expressions exist for the constitutive relations. The Drude-Born-Fedorov constitutive equations are used hereinafter.

Optical activity of the chiral medium can be explained by the direct substitution of the Drude-Born-Dedorov constitutive equations \cite{HA1,HA2}:
\begin{equation}
\left\{
\begin{aligned}
	D &=\epsilon(x)(E+\beta(x)\nabla\times E)\quad\mathrm{in}~\mathbb{R}^3, \\
	B &= \mu(x) (H+\beta(x)\nabla\times H)\quad\mathrm{in}~\mathbb{R}^3,
\end{aligned}
\right.
\end{equation}
into Maxwell's equations
\begin{equation}\label{eq:equationsection11}
\left\{
\begin{aligned}
\nabla\times E &= i\omega B \quad\mathrm{in}~\mathbb{R}^3, \\
\nabla\times H &= -i\omega D \quad\mathrm{in}~\mathbb{R}^3,
\end{aligned}
\right.
\end{equation}
which gives the constitutive relations 
\begin{equation}\label{eq:equationsection12}
\left\{
\begin{aligned}
	(1-(k(x)\beta(x))^2)D=\epsilon(x)E+i\frac{\beta(x)}{\omega}(k(x))^2H, \\
	(1-(k(x)\beta(x))^2)B=\mu(x) H-i\frac{\beta(x)}{\omega}(k(x))^2 E,
\end{aligned}	
\right.
\end{equation}
where 
$$
	k(x) = \omega\sqrt{\epsilon(x)\mu(x) }.
$$
Combining \eqref{eq:equationsection11} and \eqref{eq:equationsection12} leads to
\begin{equation}\label{eq:chiralmaxwell}
\left\{
\begin{aligned}
\nabla\times E &= (\gamma(x))^2\beta(x)E+i\omega\mu_m (\frac{\gamma(x)}{k(x)})^2 H\quad\mathrm{in}~\mathbb{R}^3,\\
\nabla\times H &= (\gamma(x))^2\beta(x)H-i\omega\epsilon(x)(\frac{\gamma(x)}{k(x)})^2 E\quad\mathrm{in}~\mathbb{R}^3,
\end{aligned}
\right.
\end{equation}
where
$$
	\gamma(x)^2 = \frac{k(x)^2}{1-(k(x)\beta(x))^2}.
$$
Let us denote
$$
k_m = \omega\sqrt{\epsilon_m\mu_m}, \quad \gamma_m = \frac{k_m^2}{1-k_m^2\beta_m^2}.
$$
Throughout this paper, we assume $k_m\beta_m<1$. Note that, $k(x)=k_m$ and $\gamma(x) = \gamma_m$ when $x$ is outside $\Omega$.

Let
\begin{equation}
\left\{
	\begin{aligned}
	\gamma_1&=\frac{\omega\sqrt{\epsilon_m\mu_m}}{1-\omega\sqrt{\epsilon_m\mu_m}\beta_m}, \quad\omega_1 = \frac{\omega}{1-\omega\sqrt{\epsilon_m\mu_m}\beta_m}, \\
	\gamma_2&=\frac{\omega\sqrt{\epsilon_m\mu_m}}{1+\omega\sqrt{\epsilon_m\mu_m}\beta_m}, \quad\omega_2 = \frac{\omega}{1+\omega\sqrt{\epsilon_m\mu_m}\beta_m}. \\
	\end{aligned}
\right.
\end{equation}

Consider an incident plane wave given by
\begin{equation}
\left\{
\begin{aligned}
E^{in}(x) &= q_1e^{i\gamma_1p_1\cdot x} + q_2e^{i\gamma_2p_2\cdot x}, \\
H^{in}(x) &= -i\sqrt\frac{\epsilon_m}{\mu_m}(q_1e^{i\gamma_1p_1\cdot x} + q_2e^{i\gamma_2 p_2\cdot x}),
\end{aligned}
\right.
\end{equation}
where
the complex vectors $p_1$, $p_2$, $q_1$, and $q_2$ satisfy the following relations:  
\begin{equation}\label{eq:incidentassumption}
\begin{aligned}
	p_1\cdot q_1=0&,\quad p_1\times q_1=-iq_1, \\
	p_2\cdot q_2=0&, \quad p_2\times q_2=iq_2.
\end{aligned}
\end{equation}
Note that, under assumption \eqref{eq:incidentassumption}, we have
\begin{equation}
\left\{
\begin{aligned}
	\nabla\times (q_1e^{i\gamma_1p_1\cdot x}) &= \gamma_1q_1e^{i\gamma_1p_1\cdot x}, \\
	\nabla\times(q_2e^{i\gamma_1p_2\cdot x}) &= -\gamma_2q_2e^{i\gamma_1p_2\cdot x}.
\end{aligned}
\right.
\end{equation}
The incident field $E^{in}$ is then a combination of a left-circularly polarized plane wave and a right-circularly polarized one, and it satisfies the homogeneous Drude-Born-Fedorov equations in $\mathbb{R}^3$.
\subsection{Radiation condition and Lippmann-Schwinger representation}

Let $E^{sc}:= E- E^{in}$ and $H^{sc}:= H-H^{in}$ be the scattered electric and magnetic fields, respectively.  In \cite{HA1}, it is established that the classical  Silver-M\"uller radiation condition, 
$$
	\big| E^{sc}(x) -\sqrt{\frac{\mu_m}{\epsilon_m}} H^{sc}(x)\times \frac{x}{|x|} \big|\leq \frac{C}{|x|^2}\quad\mathrm{for}~|x|\rightarrow +\infty,
$$
uniformly in $x/|x|$, remains valid in chiral media.  Moreover, there exists a unique solution to the scattering problem. The uniqueness follows from the Bohren decomposition of the electric and magnetic fields, 
$$
\begin{pmatrix} E \\ H \end{pmatrix} = \frac{1}{2} \begin{pmatrix} E^{(1)} \\ H^{(1)} \end{pmatrix} + \frac{1}{2} \begin{pmatrix} E^{(2)} \\ H^{(2)} \end{pmatrix},
$$
where 
$$
\begin{pmatrix} E^{(1)} \\ H^{(1)} \end{pmatrix} = \begin{pmatrix}
1 & i\sqrt\frac{\mu_m}{\epsilon_m} \\
-i\sqrt\frac{\epsilon_m}{\mu_m} & 1
\end{pmatrix} \begin{pmatrix} E \\ H \end{pmatrix} \quad \mbox{and } \quad  \begin{pmatrix} E^{(2)} \\ H^{(2)} \end{pmatrix} = \begin{pmatrix}
1 & -i\sqrt\frac{\mu_m}{\epsilon_m} \\
i\sqrt\frac{\epsilon_m}{\mu_m} & 1
\end{pmatrix} \begin{pmatrix} E \\ H \end{pmatrix}.
$$
 The existence is established by using an integral equation approach; see \cite[Theorem 5.6]{HA1}.

Let us introduce the fundamental solution to the isotropic Drude-Born-Fedorov
chiral medium.
Let $g^k$ be the outgoing fundamental solution of $\Delta + k^2$, i.e., 
$$
	g^k(x) = \frac{e^{ik |x|}}{4\pi|x|} \quad \mbox{for} \quad x\neq 0.
$$
According to \cite{HA1}, the outgoing fundamental solution $G$ of \eqref{eq:chiralmaxwell} is given by
$$
G = \frac{1}{2}\left(\mathcal{G}_1\begin{pmatrix}
1 &\ds i\sqrt\frac{\mu_m}{\epsilon_m} \\
-i\sqrt\frac{\epsilon_m}{\mu_m} & 1
\end{pmatrix} + \mathcal{G}_2\begin{pmatrix}
1 &\ds -i\sqrt\frac{\mu_m}{\epsilon_m} \\
i\sqrt\frac{\epsilon_m}{\mu_m} & 1
\end{pmatrix}\right),
$$
with
$$
\mathcal{G}_j = \frac{\gamma_j^2}{\omega_j}\begin{pmatrix}
\ds
1 + \frac{\nabla\nabla\cdot}{\gamma_j^2} &
\ds i\sqrt\frac{\mu_m}{\epsilon_m}\frac{1}{\gamma_j}\nabla\times \\
\ds
-i\sqrt\frac{\epsilon_m}{\mu_m}\frac{1}{\gamma_j}\nabla\times &
\ds 1 + \frac{\nabla\nabla\cdot}{\gamma_j^2}
\end{pmatrix}g^{\gamma_j}.
$$

In \cite{HA1}, it was proved that the following Lippmann-Schwinger representation formula for Drude-Born-Fedorov equation \eqref{eq:chiralmaxwell} holds:
\begin{equation}\label{eq:chiralgreen}
\begin{pmatrix}
E(x)-E^{in}(x) \\ H(x)-H^{in}(x)
\end{pmatrix}
= \omega\int_{\Omega}G(x-y)
\begin{pmatrix}
\tilde{\epsilon}(y) & i\omega\tilde{\tilde{\mu}}(y) \\
-i\omega\tilde{\tilde{\epsilon}}(y) & \tilde{\mu}(y)
\end{pmatrix}
\begin{pmatrix}
E(y) \\ H(y)
\end{pmatrix}
\mathrm{d}y,
\end{equation}
where 
\begin{align*}
\tilde{\epsilon}(y)& = \frac{\epsilon(y)}{\epsilon_m(1-\omega^2\epsilon(y)\mu(y) \beta^2(y))} - \frac{1}{1-\omega^2\epsilon_m\mu_m\beta_m^2}, \\
\tilde{\mu}(y) &= \frac{\mu(y)}{\mu_m(1-\omega^2\epsilon(y)\mu(y) \beta^2(y))} - \frac{1}{1-\omega^2\epsilon_m\mu_m\beta_m^2}, \\
\tilde{\tilde{\epsilon}}(y) &=  \frac{\epsilon(y) \mu(y)\beta(y)}{\mu_m(1-\omega^2\epsilon(y)\mu(y) \beta^2(y))}-\frac{\epsilon_m\beta_m}{1-\omega^2\epsilon_m\mu_m\beta_m^2}, \\
\tilde{\tilde{\mu}}(y)& =  \frac{\epsilon(y)\mu(y) \beta(y)}{\epsilon_m(1-\omega^2\epsilon(y)\mu(y) \beta^2(y))}-\frac{\mu_m\beta_m}{1-\omega^2\epsilon_m\mu_m\beta_m^2} .
\end{align*}

\section{Derivation of the dipole approximation} \label{sec2}

In this section, by the method of matched asymptotic expansions, we construct asymptotic expansions  of the scattered electromagnetic fields by the particle $\Omega$ as its characteristic size  $\delta \rightarrow 0$.  
We prove that, as $\delta \rightarrow 0$,  the particle can be approximated by a pair of  electric and magnetic dipole sources. We then show that the sources can be resonant for some negative permittivity $\epsilon_c$. It is worth emphasizing that, in a non-chiral medium, a dielectric plasmonic particle acts only as an electric dipole source \cite{pierre, vogelius}.  

\subsection{Matched asymptotic expansion}

As in \cite{khelifi}, to reveal the nature of the perturbations in the electric and magnetic fields, we introduce the local variables $\xi=x/\delta$ and set the fields $e_\delta(\xi) = E(\delta\xi)$ and $h_\delta(\xi) = H(\delta\xi)$. We expect that $E(x)$ and $H(x)$ will differ appreciably from $E^{in}(x)$ and  $H^{in}(x)$ for $x$ close to $0$, but they will differ little from $E^{in}(x)$ and $H^{in}(x)$ for $x$ far from $0$. Therefore, in the spirit of matched asymptotic expansions, we shall represent each of the fields $E$ and $H$ by two different expansions, an {inner expansion} for $x$ near $0$, and an {outer expansion} for $x$ far from $0$. The outer expansions must begin with $E^{in}$ and $H^{in}$, so we write:
\begin{equation}\label{eq:outerexpansion}
\begin{aligned}
	E(x)&=E^{in}(x) + \delta^{\alpha_1}E_1(x) + \delta^{\alpha_2}E_2(x) + \cdots, \quad\mathrm{for}~|x|\gg O(\delta) , \\
	H(x)&=H^{in}(x) + \delta^{\alpha_1}H_1(x) + \delta^{\alpha_2}H_2(x) + \cdots, \quad\mathrm{for}~|x|\gg O(\delta),
\end{aligned}
\end{equation}
where $0<\alpha_1<\alpha_2<\cdots$, and $(E_1, H_1)$, $(E_2, H_2)$, $\cdots$ are to be found. Inserting this series into equation \eqref{eq:chiralmaxwell} and observing that
$$
	\epsilon(\frac{x}{\delta})\equiv \epsilon_m,\quad \mu(\frac{x}{\delta})\equiv \mu_m
$$
for $|x|\gg O(\delta)$, we find that the outer coefficients $(E_j, H_j)$, $j=1,2,\cdots$ are  solutions to
$$
\left\{
\begin{aligned}
\nabla\times E_j &= \gamma_m^2\beta_m E_j +i\omega\mu_m (\frac{\gamma_m}{k_m})^2 H_j,\\
\nabla\times H_j &= \gamma_m^2\beta_m H_j -i\omega\epsilon_m (\frac{\gamma_m}{k_m})^2 E_j
\end{aligned}
\right.
$$
for $|x| \gg O(\delta)$. Moreover, all $(E_j, H_j)$ satisfy the Silver-M\"uller radiation condition. 


We write the inner expansion as
\begin{equation}\label{eq:innerexpansion}
\begin{aligned}
E(\delta\xi) &= e_\delta(\xi) = e_0(\xi) + \delta e_1(\xi) + \cdots \quad\mathrm{for}~|\xi|=O(1),\\
H(\delta\xi) &= h_\delta(\xi) = h_0(\xi) + \delta h_1(\xi) + \cdots \quad\mathrm{for}~|\xi|=O(1),
\end{aligned}
\end{equation}
where $(e_0,h_0)$, $(e_1,h_1)$, $(e_2,h_2)$, $\cdots$ are to be found. The functions $e_j(\xi)$ and $h_j(\xi)$ are defined everywhere in $\mathbb{R}^3$. 

In order to determine $E_j(x), H_j(x), e_j(\xi)$, and $h_j(\xi)$, we  need to equate the inner and the outer expansions in some overlap domain within which the stretched variable $\xi$ is large and $x$ is small. In this domain the matching conditions are
\begin{equation}\label{eq:matchingconditionold}
\begin{aligned}
E^{in}(x) + \delta^{\alpha_1} E_1(x) + \cdots \sim e_0(\xi) + \delta e_1(\xi) + \cdots, \\
H^{in}(x) + \delta^{\alpha_1} H_1(x) + \cdots \sim h_0(\xi) + \delta h_1(\xi) + \cdots.
\end{aligned}
\end{equation}
These matching conditions will be made more precise later on.

Inserting the inner expansions into the Lippmann-Schwinger representation formula \eqref{eq:chiralgreen}, we arrive at
\begin{equation}
\begin{pmatrix}
 \delta^{\alpha_1}E_1(x) + \cdots \\ \delta^{\alpha_1}H_1(x) + \cdots
\end{pmatrix}
= \omega \delta^3 \int_{\Omega}G(x- \delta \xi)
\begin{pmatrix}
	\tilde{\epsilon}(\delta \xi) & i\omega\tilde{\tilde{\mu}}(\delta \xi) \\
	-i\omega\tilde{\tilde{\epsilon}}(\delta \xi) & \tilde{\mu}(\delta \xi)
\end{pmatrix}
\begin{pmatrix}
	e_0(\xi)+\delta e_1(\xi) + \cdots \\ h_0(\xi)+\delta h_1(\xi) + \cdots
\end{pmatrix}
\mathrm{d}\xi
\end{equation} 
for $|x| \gg O(\delta)$. 
Therefore, the matching conditions (\ref{eq:matchingconditionold}) imply that
$$
	\alpha_i = 2+i,\quad i=1,2,\cdots.
$$

If we substitute the inner expansion \eqref{eq:innerexpansion} into \eqref{eq:chiralmaxwell} and equate coefficients of $\delta^{-1}$, then we get the following equations:
\begin{equation}\label{eq:matching}
\left\{
	\begin{aligned}
		&\nabla \times e_0 = 0	\quad\mathrm{in}~\mathbb{R}^3, \\
		&\nabla\cdot\Big(\epsilon(1+\gamma^2\beta^2)e_0+i\omega \epsilon\mu\beta \Big(\frac{\gamma}{k}\Big)^2h_0\Big) = 0\quad\mathrm{in}~\mathbb{R}^3, \\
		&\nabla\times h_0 = 0\quad\mathrm{in}~\mathbb{R}^3, \\
		&\nabla\cdot\Big(\mu(1+\gamma^2\beta^2)h_0-i\omega \epsilon\mu\beta \Big(\frac{\gamma}{k}\Big)^2e_0\Big) = 0\quad\mathrm{in}~\mathbb{R}^3,\\
		&e_0(\xi)\rightarrow E^{in}(0)\quad\mathrm{as}~|\xi|\rightarrow +\infty, \\
		&h_0(\xi)\rightarrow H^{in}(0)\quad\mathrm{as}~|\xi|\rightarrow +\infty.
	\end{aligned}
\right.
\end{equation}
Here, the derivatives are taken with respect to $\xi$.

Since the curl of $e_0(\xi)$ and $h_0(\xi)$ are both zero, there exists scalar functions $V(\xi)$ and $W(\xi)$  satisfying
$$
	e_0(\xi) = \nabla V(\xi),\qquad h_0(\xi) = \nabla W(\xi), \quad \xi \in \mathbb{R}^3.
$$
Then \eqref{eq:matching} becomes
\begin{equation}\label{eq:matching1}
\left\{
\begin{aligned}
&\nabla\cdot\Big(\epsilon(1+\gamma^2\beta^2)\nabla V+i\omega \epsilon\mu\beta \Big(\frac{\gamma}{k}\Big)^2\nabla W\Big) = 0\quad\mathrm{in}~\mathbb{R}^3, \\
&\nabla\cdot\Big(\mu(1+\gamma^2\beta^2)\nabla W-i\omega \epsilon\mu\beta \Big(\frac{\gamma}{k}\Big)^2\nabla V\Big) = 0\quad\mathrm{in}~\mathbb{R}^3,\\
& V(\xi)-E^{in}(0)\cdot \xi\rightarrow 0\quad\mathrm{as}~|\xi|\rightarrow +\infty, \\
& W(\xi)-H^{in}(0)\cdot \xi\rightarrow 0\quad\mathrm{as}~|\xi|\rightarrow +\infty.
\end{aligned}
\right.
\end{equation}
The first two equations in \eqref{eq:matching1} imply that the normal components of 
$$
\epsilon(1+\gamma^2\beta^2)\nabla V+i\omega \epsilon\mu\beta \Big(\frac{\gamma}{k}\Big)^2\nabla W,
$$
and
$$
\mu(1+\gamma^2\beta^2)\nabla W-i\omega \epsilon\mu\beta \Big(\frac{\gamma}{k}\Big)^2\nabla V,
$$
are continuous on $\partial B$. Therefore, \eqref{eq:matching1} is equivalent to
\begin{equation}\label{eq:equationvw}
\left\{
\begin{aligned}
	&\Delta V = \Delta W = 0	\quad\mathrm{in}~B\cup(\mathbb{R}^3 \backslash \bar{B}), \\
	& V|^- =  V|^+,\quad  W|^- =  W|^+\quad\mathrm{on}~\partial B, \\
	&\epsilon_c\left.\frac{\partial V}{\partial\nu}\right|^- = 
\epsilon_m (1 + \gamma_m^2 \beta_m^2) 	\left.\frac{\partial V}{\partial\nu}\right|^+
+ i \omega \epsilon_m \mu_m \beta_m \left(\frac{\gamma_m}{k_m}\right)^2 
\left.\frac{\partial W}{\partial\nu}\right|^+\quad\mathrm{on}~\partial B, \\
	&\mu_c\left.\frac{\partial W}{\partial\nu}\right|^- = \mu_m (1 + \gamma_m^2 \beta_m^2) 	\left.\frac{\partial W}{\partial\nu}\right|^+ - i \omega \epsilon_m \mu_m \beta_m 
	\left(\frac{\gamma_m}{k_m}\right)^2 
\left.\frac{\partial V}{\partial\nu}\right|^+
\quad\mathrm{on}~\partial B,\\
	&( V -  V^{in})(\xi)\rightarrow 0\quad\mathrm{as}~|\xi|\rightarrow\infty, \\
	&( W -  W^{in})(\xi)\rightarrow 0\quad\mathrm{as}~|\xi|\rightarrow\infty.
\end{aligned}
\right.
\end{equation}
Here, $\frac{\partial}{\partial \nu}$ denotes the normal derivative on $\partial B$, and the subscripts $+$ and $-$ indicate the limits from outside and inside $B$, respectively.

\subsection{Small volume expansion and its resonant behavior}

In this subsection, we formally derive a small volume expansion of $(E,H)$ by solving $(e_0,h_0)$ and using the Lippmann-Schwinger equation. We emphasize that the expansion can be rigorously proved.
To solve $e_0=\nabla V$ and $h_0=\nabla W$, we make use of boundary layer potentials. We also discuss the resonant behavior of the expansion due to the negative permittivity of the plasmonic particle.

We define the single-layer potential $\mathcal{S}_B$ as
\begin{equation*}
\mathcal{S}_{B}[\psi](x):= - \frac{1}{4 \pi}\int_{\partial B} \frac{1}{ |x-y|}\psi(y) d\sigma(y),~~~x\in \mathbb{R}^{3},
\end{equation*}
for  $\psi\in H^{-1/2}(\partial B)$.
We also define the Neumann-Poincar\'e operator  $\mathcal{K}_{B}^{*}$ by
\begin{equation*}
\mathcal{K}_{B}^{*}[\psi](x) =\int_{\partial B}\frac{(x -y)\cdot \nu(x)}{4 \pi |x-y|^3} \psi(y)d\sigma(y),~~~x\in \partial B
\end{equation*}
with $\nu(x)$ being the outward normal at $x \in \partial B$.
It is well-known that the following jump relation holds:
\begin{equation} \label{jump}
	\left.\frac{\partial \mathcal{S}_B}{\partial\nu}\right|_{\partial B}^{\pm}[\psi] = (\pm\frac{1}{2} + \mathcal{K}_B^*)[\psi].
\end{equation}
 
The functions $V$ and $W$ can be represented by using the single layer potential as
\begin{equation}\label{eq:vwpsiepsih}
\begin{aligned}
V &= V^{in}+ \mathcal{S}_B[\psi_E], \\
W &= W^{in}+ \mathcal{S}_B[\psi_H].
\end{aligned}
\end{equation}
Using the transmission conditions on $\partial B$ in \eqref{eq:equationvw} and the jump relation \eqref{jump}, it can be shown that the pair $(\psi_E,\psi_H)$ is the solution to the boundary integral equation
\begin{equation}\label{eq:psiEpsiH}
	A(\epsilon_c)
	\begin{pmatrix}
	\psi_E \\ \psi_H
	\end{pmatrix}
	=\begin{pmatrix}
	f_E \\ f_H
	\end{pmatrix},
\end{equation}
where 
$$
\mathcal{A} =\begin{pmatrix}
\lambda_\epsilon I- \mathcal{K}_B^* & i\omega d_\epsilon(\frac{1}{2}+ \mathcal{K}_B^*) \\
-i\omega d_\mu(\frac{1}{2}+ \mathcal{K}_B^*) & \lambda_\mu I- \mathcal{K}_B^*
\end{pmatrix},
$$
and
\begin{align}
\begin{pmatrix}
f_E \\ f_H
\end{pmatrix} &:= \begin{pmatrix}
1 & -i\omega d_\epsilon \\
i\omega d_\mu & 1
\end{pmatrix}
\begin{pmatrix}
E^{in}(0)\cdot\nu \\
H^{in}(0)\cdot\nu
\end{pmatrix}\bigg|_{\partial B}. 
\label{eq:deffefh}
\end{align}
Here, the parameters $\lambda_\epsilon,\lambda_\mu,d_\epsilon$ and $d_\mu$ are given by
\begin{align*}
\lambda_\epsilon &= \frac{\epsilon_c+\epsilon_m(1+{\gamma_m^2}{\beta_m^2})}{2(\epsilon_c-\epsilon_m(1+{\gamma_m^2}{\beta_m^2}))}, \\
\lambda_\mu &= \frac{\mu_m + \mu_m(1+{\gamma_m^2}{\beta_m^2})}{2(\mu_m-\mu_m(1+{\gamma_m^2}{\beta_m^2}))}, \\
d_\epsilon &= \frac{\epsilon_m\mu_m \beta_m(\gamma_m/k_m)^2}{\epsilon_c-\epsilon_m(1+{\gamma_m^2}{\beta_m^2})}, \\
d_\mu &=\frac{\epsilon_m \mu_m\beta_m(\gamma_m/k_m)^2}{\mu_m-\mu_m(1+{\gamma_m^2}{\beta_m^2})}.
\end{align*}

It is known that the operator $\mathcal{K}_B^*$ can be symmetrized using Calder\'on identity and hence becomes self-adjoint  \cite{HA3,hyeonbae1}. Let $H_0^{-1/2}(\p B)$ be the subspace of $H^{-1/2}(\p B)$ with zero mean value. Let $\mathcal{H}^*(\partial B)$ be the space $H_0^{-1/2}(\partial B)$ equipped with inner product
$$
\langle \varphi_1,\varphi_2\rangle_* = -\langle \mathcal{S}_B[\varphi_2],\varphi_1\rangle_{\frac{1}{2},-\frac{1}{2}},
$$
and $(\lambda_j,\phi_j)$, $j=1,2,\cdots$ be the pair of eigenvalue and normalized eigenfunction of $\mathcal{K}_B^*$ in $\mathcal{H}^*(\p B)$. For any $\varphi\in H^*(\partial B)$, the following spectral representation formula holds:
$$
	\mathcal{K}_B^*[\varphi] = \sum\limits_{j=1}^\infty \lambda_j\phi_j {\langle \phi_j,\varphi\rangle_*}.
$$
It is then easy to see that the operator $\mathcal{A}$ has a block matrix structure. Indeed,  we have 
\begin{align}
\mathcal{A} \begin{pmatrix}a \phi_{n} \\ b \phi_{n} \end{pmatrix}&=
\begin{pmatrix}
\lambda_\epsilon(\epsilon_c)  - \lambda_n & i\omega d_\epsilon(\epsilon_c)(\frac{1}{2} + \lambda_n) \\
-i\omega d_\mu(\frac{1}{2}+\lambda_n) & \lambda_\mu -\lambda_n
\end{pmatrix}
\begin{pmatrix}
a\phi_n \\ b\phi_{n}
\end{pmatrix}
\nonumber
\\&:=  A_n\begin{pmatrix}a\phi_n \\ b\phi_{n}\end{pmatrix}.
\label{eq:defan}
\end{align}

Since $\{\phi_j\}_{j=1}^\infty$ forms a complete and orthonormal basis, we can write
\begin{align*}
	&\psi_E=\sum\limits_{n=1}^\infty \psi_E^n\phi_n,\quad \psi_H=\sum\limits_{n=1}^\infty \psi_H^n\phi_n, \\
	&f_E = \sum\limits_{n=1}^\infty f_E^n\phi_n,\quad f_H=\sum\limits_{n=1}^\infty f_H^n\phi_n,
\end{align*}
with the coefficients
\begin{align*}
	&\psi_E^n = \langle\psi_E,\phi_n\rangle_*,\quad \psi_H^n = \langle\psi_H,\phi_n\rangle_*, \\
	&f_E^n = \langle f_E,\phi_n\rangle_*, \quad f_H^n = \langle f_H,\phi_n\rangle_*.
\end{align*}
Then, from \eqref{eq:defan}, we have
\begin{equation}\label{eq:phiephih1}
	\begin{pmatrix}
	\psi_E^n \\ \psi_H^n
	\end{pmatrix} ={A}_n^{-1}\begin{pmatrix}
	f_E^n \\ f_H^n
	\end{pmatrix}={A}_n^{-1}\begin{pmatrix}
	1 & -i\omega d_\epsilon \\
	i\omega d_\mu & 1
	\end{pmatrix}
	\begin{pmatrix}
	\langle E^{in}(0)\cdot\nu, \phi_n\rangle_* \\
	\langle H^{in}(0)\cdot\nu, \phi_n\rangle_*
	\end{pmatrix}\bigg|_{\partial B}.
\end{equation}
Therefore, we obtain the solution $V$ and $W$ in terms of the eigenvalue $\lambda_n$ and eigenfunctions $\phi_n$ of the Neumann-Poincar\'e operator $\mathcal{K}^*_B$.

Now, we derive a small volume expansion of $(E,H)$ using the Lippmann-Schwinger equation. 
Since $\Omega=\delta B$ and $\delta\ll 1$, by Taylor's expansion,  \eqref{eq:chiralgreen} leads to 
\begin{equation}\label{eq:lippmanschwingera_temp}
\begin{aligned}
\begin{pmatrix}
E-E^{in} \\ H-H^{in}
\end{pmatrix}(x)
&= \omega\int_{\Omega}(G(x)+O(\delta))\begin{pmatrix}
\tilde{\epsilon}(0) & i\omega\tilde{\tilde{\mu}}(0) \\
-i\omega\tilde{\tilde{\epsilon}}(0) & \tilde{\mu}(0)
\end{pmatrix}\begin{pmatrix}
E(y) \\ H(y)
\end{pmatrix}
\mathrm{d}y \\
&\approx \omega G(x) K_0\int_\Omega\begin{pmatrix}
E(y) \\ H(y)
\end{pmatrix}\mathrm{d}y \\
&\approx \omega G(x)K_0 \delta^3\int_B\begin{pmatrix}
e_0 \\ h_0 \end{pmatrix},
\end{aligned}
\end{equation}
where  $K_0$ is defined by 
$$
K_0= \begin{pmatrix}
 \ds \frac{\epsilon_c}{\epsilon_m} - \frac{1}{1-\omega^2\epsilon_m\mu_m\beta_m^2}& 
 \ds -i\omega\frac{\mu_m\beta_m}{1-\omega^2\epsilon_m\mu_m\beta_m^2} \\
\ds  i\omega \frac{\epsilon_m\beta_m}{1-\omega^2\epsilon_m\mu_m\beta_m^2} & 
 \ds 1-  \frac{1}{1-\omega^2\epsilon_m\mu_m\beta_m^2}
\end{pmatrix}.
$$
Since
\begin{equation}
\begin{pmatrix}
e_0 \\ h_0
\end{pmatrix}=
\begin{pmatrix}
E^{in}(0) \\ H^{in}(0)
\end{pmatrix}
+ \begin{pmatrix}
\nabla \mathcal{S}_B[\psi_E] \\
\nabla \mathcal{S}_B[\psi_H]
\end{pmatrix}, 
\end{equation}
we get
\begin{equation}\label{eq:lippmanschwingera}
\begin{aligned}
\begin{pmatrix}
E-E^{in} \\ H-H^{in}
\end{pmatrix}(x)
&\approx 
\omega G(x)K_0(\epsilon_c) \delta^3\int_B\begin{pmatrix}
E^{in}(0) + \nabla \mathcal{S}_B[\psi_E] \\ H^{in}(0)+\nabla \mathcal{S}_B[\psi_H] \end{pmatrix}
.
\end{aligned}
\end{equation}

Hence, we need to analyze the integral
$$
\int_B\begin{pmatrix}
\nabla \mathcal{S}_B[\psi_E] \\ \nabla \mathcal{S}_B[\psi_H]
\end{pmatrix}.
$$
Let us define
\begin{equation}\label{def:mn}
\begin{pmatrix}
	M_n^{EE} & M_n^{EH} \\
	M_n^{HE} & M_n^{HH}
\end{pmatrix}=(-1){A}_n^{-1}
\begin{pmatrix}
	1 & -i \omega d_\epsilon \\
	i\omega d_\mu & 1
\end{pmatrix}.
\end{equation}
Then \eqref{eq:phiephih1} can be written as
\begin{align*}
\psi_E^n &= -M_n^{EE}\langle E^{in}(0)\cdot\nu,\phi_n\rangle_* - M_n^{EH}\langle H^{in}(0)\cdot\nu, \phi_n\rangle_*, \\
\psi_H^n &= -M_n^{HE}\langle E^{in}(0)\cdot\nu,\phi_n\rangle_* - M_n^{HH}\langle H^{in}(0)\cdot\nu, \phi_n\rangle_*.
\end{align*}
For convenience of notation, here we slightly generalize the definition of the inner product $\langle \varphi_1,\varphi_2\rangle_*$ to the case when $\varphi_1\in(\mathcal{H}^*(\partial B))^k,k\in\mathbb{N}$ by
$$
	\langle\varphi_1,\varphi_2\rangle_* = \begin{pmatrix}
	\langle\varphi_1^1,\varphi_2\rangle_* \\
	\langle\varphi_1^2,\varphi_2\rangle_* \\
	\vdots \\
	\langle\varphi_1^k,\varphi_2\rangle_* \\
	\end{pmatrix}.
$$
With this notation, we get
\begin{align*}
\psi^n_E &= -M_n^{EE}\langle \nu,\phi_n \rangle^\top_*E^{in}(0) - M_n^{EH}\langle \nu,\phi_n\rangle^\top_*H^{in}(0), \\
\psi^n_H &= -M_n^{HE}\langle \nu,\phi_n \rangle^\top_*E^{in}(0) - M_n^{HH}\langle \nu,\phi_n\rangle^\top_*H^{in}(0),
\end{align*}
where the superscript $\top$ denotes the Hermitian conjugate. By using the integration by parts, it follows that
\begin{equation}\label{eq:psie1}
\begin{aligned}
\int_B \nabla \mathcal{S}_B[\psi_E]\mathrm{d}x &= \sum\limits_{n=1}^\infty \psi_E^n\int_B\nabla \mathcal{S}_B[\phi_n]\mathrm{d}x \\
&= \sum\limits_{n=1}^\infty(-1)\psi_E^n\langle \nu,\phi_n\rangle_* \\
&= \sum\limits_{n=1}^\infty M_n^{EE}\langle\nu,\phi_n\rangle_*\langle\nu,\phi_n\rangle_*^\top E^{in}(0) + \sum\limits_{n=1}^\infty M_n^{EH}\langle\nu,\phi_n\rangle_*\langle\nu,\phi_n\rangle_*^\top H^{in}(0),
\end{aligned}
\end{equation}
and similarly for $\psi_H$, 
\begin{equation}\label{eq:psih1}
\begin{aligned}
	\int_B \nabla \mathcal{S}_B[\psi_H]\mathrm{d}x= \sum\limits_{n=1}^\infty M_n^{HE}\langle\nu,\phi_n\rangle_*\langle\nu,\phi_n\rangle_*^\top E^{in}(0) + \sum\limits_{n=1}^\infty M_n^{HH}\langle\nu,\phi_n\rangle_*\langle\nu,\phi_n\rangle_*^\top H^{in}(0).
\end{aligned}
\end{equation}
Now, if we define the polarization tensor by
\begin{equation}\label{eq:polarizationdefine}
	M=M(\epsilon_c,B)=\begin{pmatrix}
	M^{EE} & M^{EH} \\
	M^{HE} & M^{HH}
	\end{pmatrix},
\end{equation}
where
\begin{equation}\label{eq:eeehhehh}
\begin{aligned}
M^{EE} &= \sum\limits_{n=1}^\infty M_n^{EE}\langle\nu,\phi_n\rangle_*\langle\nu,\phi_n\rangle_*^\top, \\
M^{EH} &= \sum\limits_{n=1}^\infty M_n^{EH}\langle\nu,\phi_n\rangle_*\langle\nu,\phi_n\rangle_*^\top, \\
M^{HE} &= \sum\limits_{n=1}^\infty M_n^{HE}\langle\nu,\phi_n\rangle_*\langle\nu,\phi_n\rangle_*^\top, \\
M^{HH} &= \sum\limits_{n=1}^\infty M_n^{HH}\langle\nu,\phi_n\rangle_*\langle\nu,\phi_n\rangle_*^\top,
\end{aligned}
\end{equation}
then \eqref{eq:psie1} and \eqref{eq:psih1} can be rewritten  as
\begin{equation}\label{eq:equation3}
\int_B \begin{pmatrix}
\nabla \mathcal{S}_B[\psi_E] \\ \nabla \mathcal{S}_B[\psi_H]
\end{pmatrix}\mathrm{d}x = M(\epsilon_c,B)
\begin{pmatrix}
	E^{in}(0) \\
	H^{in}(0)
\end{pmatrix}.
\end{equation}

Finally, from \eqref{eq:lippmanschwingera} and \eqref{eq:equation3}, we obtain a small volume expansion
\begin{equation}\label{eq:lippmanschwinger}
\begin{aligned}
\begin{pmatrix}
E-E^{in} \\ H-H^{in}
\end{pmatrix}(x)
&\approx 
\omega G(x)K_0(\epsilon_c) \delta^3 \widetilde{M}(\epsilon_c,B) \begin{pmatrix} 
E^{in}(0)  \\ H^{in}(0) \end{pmatrix}
,\quad |x| \gg O(\delta),
\end{aligned}
\end{equation}
where $\widetilde{M}:= |B| I + M$.

Now, let us discuss a resonant behavior of the polarization tensor $M$.
Straightforward computation shows that 
$$
 \det A_n (\epsilon_c)= (-1)\frac{(1/2-\lambda_n)(1- k_m^2\beta_m^2 (1/2-\lambda_n)) ({1-k_m^2\beta_m^2}) }{k_m^2\beta_m^2 } \frac{\epsilon_c-\epsilon_{c,n}^*}{({1-k_m^2\beta_m^2})\epsilon_c -  \epsilon_m} , 
$$
where
$$
\epsilon_{c,n}^* = -\epsilon_m\frac{\ds 1/2+\lambda_n}{\ds 1/2 -\lambda_n }\big({1-k_m^2 \beta_m^2 \big(1/2-\lambda_n\big)}\big)^{-1}.
$$
It is clear that $\epsilon_{c,n}^*<0$. Therefore, when the particle is plasmonic, i.e., $\mathrm{Re}\{\epsilon_c\}<0$, the polarization tensor can be very large if $\epsilon_c$ is close to $\epsilon_{c,n}^*<0$ for some $n$.

Regarding the permittivity $\epsilon_c$, We make the following assumptions.
\begin{asump}\label{asump_epc}
Suppose that
\begin{itemize}
\item[(i)] There exists $n\in\mathbb{N}$ such that $\langle\nu,\phi_n\rangle_* \neq 0$;
\item[(ii)] The permittivity $\epsilon_c$ of the plasmonic particle  is close to $\epsilon_{c,n}^*<0$. \end{itemize}
\end{asump}
It is worth mentioning that, when $\Omega$ is a ball, then the above assumption is satisfied with $n=1$ (see Appendix \ref{app_ball}).
Under the above assumption, $A_{n}(\epsilon_c)^{-1}$ is nearly singular. In view of (\ref{def:mn}), it follows that
$$
\widetilde{M}(\epsilon_c,B) = M_n(\epsilon_c,B) + O(1) = O\Big(\frac{1}{\det A_n(\epsilon_c)}\Big).
$$
Therefore, from \eqref{eq:lippmanschwinger}, we arrive at
\begin{equation} \label{formal_SVE}
\begin{aligned}
\begin{pmatrix}
E(x)-E^{in}(x) \\
H(x)-H^{in}(x)
\end{pmatrix} &\approx \delta^3  \omega G(x)K_0 {M}_n(\epsilon_c,B) \langle\nu,\phi_n\rangle_*\langle\nu,\phi_n\rangle_*^\top\begin{pmatrix}
E^{in}(0) \\ H^{in}(0)
\end{pmatrix}.
\end{aligned}
\end{equation}

\begin{theorem} \label{thmsmall} For $|x| \gg O(\delta)$, the following asymptotic expansion holds:
\begin{align} 
\begin{pmatrix}
E(x)-E^{in}(x) \\
H(x)-H^{in}(x)
\end{pmatrix} &= \delta^3  \omega G(x)K_0 {M}_n(\epsilon_c, B) \langle\nu,\phi_n\rangle_*\langle\nu,\phi_n\rangle_*^\top \begin{pmatrix}
E^{in}(0) \\ H^{in}(0)
\end{pmatrix}
  +O(\frac{\delta^4}{|\mathrm{det}A_n|}).
 \label{asympform}
\end{align}
\end{theorem}

\begin{rmk}
Using the layer potential formulation (\ref{eq**}) in Appendix \ref{subsec:layer_chiral}, Theorem \ref{thmsmall} can be proved rigorusly by applying  to (\ref{eq**}) essentially the same method  as the one in \cite{pierre}. 

\end{rmk}

\begin{rmk}
Note that $K_0=0$ if $\epsilon_c=\epsilon_m$. Moreover, from the definitions of $M^{EE}_n, M^{EH}_n, M^{HE}_n, M^{HH}_n$ and $K_0$, we can conclude that ${M}_n$ and $K_0$ are independent of the position of particles. Formula \eqref{asympform} shows that the scattered  wave from a single particle behaves similarly to a pair of resonant electric and magnetic dipole sources in the far-field. These dipole sources resonate at the set of permittivities $\epsilon_c$ satisfying $\mathrm{det}A_n(\epsilon_c)\approx 0$.  
\end{rmk}

\section{Effective medium theory and double-negative materials} \label{sec_EMT}

We are now ready to study when we could reach double-negative mode with  multiple dilute nanoparticles embedded into a chiral medium. We first derive an effective medium theory for a large number nanoparticles embedded in a chiral medium. Some conditions on the volume fraction and distribution of the nanoparticles are required. Then we prove that both the effective electric permittivity and effective magnetic permeability can be negative near the resonant frequencies. 

Let ${\Omega}$ be a bounded smooth domain. We consider a collection of small identical plasmonic particles $\{\Omega_j^N\}_{j=1}^{N^3}$ with size of order $\delta$. Each particle can be represented by $\Omega_j^N = \delta B + z_j^N$, where $z_j^N$ is the center location of $\Omega_j^N$. Let 
$$
{\Omega}^N = {\Omega} \setminus \overline{\cup_{j=1}^{N^3} \Omega_j^N}.
$$

We assume that the following assumptions on the distribution  of the plasmonic particles over the domain ${\Omega}$ hold. 
\begin{asump}  \label{assumptionV1} There exists a smooth function $\widetilde{V}$ such that for arbitrary smooth functions $f$ and $g$, 
\begin{equation} \label{assumptionV2} \begin{array}{l}
\ds \max_{1\leq j\leq N^3} \bigg| \frac{1}{N^3} \sum_{i\neq j} G(z_i^N-z_j^N)   \left(\begin{array}{l} f(z_i^N) \\ g(z_i^N) \end{array}\right) - \int_{{\Omega}} G(z- z_j^N) \widetilde{V}(z)  \left(\begin{array}{l} f(z) \\ g(z) \end{array}\right)  dz \bigg| \\
\nm 
\qquad \ds 0 \mbox{ as} N\rightarrow +\infty. 
\end{array} 
\end{equation}
\end{asump}

The scattering problem of electromagnetic waves in a chiral media by the system of plasmonic particles can be modeled as 
\begin{equation}\label{eq:chiralmaxwell_N}
\left\{
\begin{aligned}
&\nabla\times E^N = i\omega\mu_m  H^N \quad\mathrm{in}~\cup_{j=1}^{N^3} \Omega_j^N,\\
&\nabla\times H^N = -i\omega\epsilon_c  E^N \quad\mathrm{in}~\cup_{j=1}^{N^3} \Omega_j^N,
\\
&\nabla\times E^N = \gamma_m^2\beta_m E^N+i\omega\mu_m \left(\frac{\gamma_m}{k_m}\right)^2 H^N \quad\mathrm{in}~\mathbb{R}^3\setminus \overline{\cup_{j=1}^{N^3} \Omega_j^N},\\
&\nabla\times H^N = \gamma_m^2\beta_m H^N-i\omega\epsilon_m\left(\frac{\gamma_m}{k_m}\right)^2 E^N \quad\mathrm{in}~\mathbb{R}^3\setminus \overline{\cup_{j=1}^{N^3} \Omega_j^N},\\
&E^N \times \nu |_- = E^N \times \nu |_+ \quad \mbox{on }\ {\partial \Omega_j^N},\  j=1,...,N^3,
\\
&H^N \times \nu |_- = H^N \times \nu |_+ \quad \mbox{on }\ {\partial \Omega_j^N},\  j=1,...,N^3.
\end{aligned}
\right.
\end{equation}
Moreover, the pair $(E^N-E^{in}, H^N-H^{in})$ satisfies the Silver-M\"uller radiation condition.

Then, using  the layer potential formulation in Appendix \ref{subsec:layer_chiral}, the solution $(E^N,H^N)$ can be represented as
$$
E^{N} =
E^{in} + \mathcal{Q}^E_{\Omega^N}\begin{bmatrix}\varphi^N \\ \psi^N\end{bmatrix},
\quad
H^{N} = H^{in} + \mathcal{Q}^H_{\Omega^N}\begin{bmatrix}\varphi^N \\ \psi^N\end{bmatrix}, 
$$
where $(\varphi^N,\psi^N)$ is the solution to
$$
(\mathcal{J}_{\Omega^N} + \mathcal{C}_{\Omega^N}) 
\begin{bmatrix}
\varphi^N
\\
\psi^N
\end{bmatrix}
 = \begin{bmatrix}
 \nu\times E^{in}|_{\partial \Omega^N}
 \\
 \nu\times H^{in}|_{\partial \Omega^N}
 \end{bmatrix}.
$$
Here, we have used the notations
\begin{align*}
&\varphi^N = (\varphi_1^N,...,\varphi_{N^3}^N ),
\\
&\psi^N = (\psi_1^N,...,\psi_{N^3}^N ),
\\
\end{align*}
and
\begin{align*}
\mathcal{F}_{\Omega^N} &\begin{bmatrix}\varphi^N \\ \psi^N\end{bmatrix} = \sum_{j=1}^{N^3}
\mathcal{F}_{\Omega^N_j}\begin{bmatrix}\varphi^N_j \\ \psi^N_j\end{bmatrix}
\end{align*}
for $\mathcal{F} = \mathcal{Q}^E, \mathcal{Q}^H, \mathcal{J}$ and  $\mathcal{C}$.

We further assume that all the particles are aligned in a dilute manner.
\begin{asump} \label{dilute2}
 There exists $\Lambda>0$ and $a>0$, such that 
$$
	\delta = \Lambda^{1/3}N^{-1-a}.
$$
\end{asump}
In this case, the total volume of all plasmonic particles is of order $O(\delta^3N^3) = O(N^{-3a})$, which converges to $0$ as  $N\rightarrow\infty$.

One of the most important reason that we align particles in dilute way is that, we can approximate the scattered field $(E-E^{in}, H-H^{in})$ by the sum of fields generated by individual 
 particles, i.e.,  the interaction between scattering fields from different particles is  negligible.

For simplicity, we also assume that the particle is symmetric so that the tensor $M_n(\epsilon_c,R_{\theta_j}B)\langle\nu,\phi_n\rangle_*\langle\nu,\phi_n\rangle_*^\top$ is proportional to the identity matrix $I$. More precisely, we assume that
\begin{equation}\label{tensor_cn}
\langle\nu,\phi_n\rangle_*\langle\nu,\phi_n\rangle_*^\top = c_n I,
\end{equation}
for some $c_n>0$.
A more general case can be considered in the same way by assuming that the particles are randomly oriented and using averaging with respect to the orientation of the particle (see Remark \ref{rmk:random}).
When $B$ is a unit ball, one can check that $c_1=\frac{4\pi}{27}$  (see Appendix \ref{app_ball}).

For $1\leq j\leq N^3$, we define  
$$
E^{in,N}_j =
E^{in} + \sum_{i\neq j} \mathcal{P}_{\Omega_i^N}\begin{bmatrix}\varphi_i^N \\ \psi_i^N\end{bmatrix},
\quad
H^{in,N}_j = H^{in} + \sum_{i\neq j}\mathcal{Q}_{\Omega^N}\begin{bmatrix}\varphi_i^N \\ \psi_i^N\end{bmatrix}.
$$
For each $1\leq j\leq N^3$ and $x\in {\Omega}^N$, 
\begin{align}
E^N(x) &= E^{in}(x) + \sum_{i=1}^{N^3} E^{sc,N}_i(x)  \nonumber
\\
&= E^{in,N}_j(x) +  E^{sc,N}_j(x).\label{EN_decomp}
\end{align}

\begin{prop}\label{prop_pointsource}
For each $1\leq j\leq N^3$ and $x\in {\Omega}^N$,
$$
\begin{pmatrix}
E^{sc,N}_j(x)
\\
H^{sc,N}_j(x)
\end{pmatrix}
=   
 \delta^3 \omega G(x-z_j^N)K_0{M}_n(\epsilon_c,B) c_n
\begin{pmatrix}
E^{in,N}(z_j^N) \\ H^{in,N}(z_j^N)
\end{pmatrix}+ (\frac{N^3 \delta^4}{ |\mathrm{det} A_n|}).
$$
\end{prop}

\subsection{Derivation of the homogenized equation and analysis of effective  parameters}

Let us assume that \begin{equation} \ds N^{-3a} |\mathrm{det} A_n|^{-1} = O(1),  \end{equation}
and define
$$
\begin{pmatrix}
\tilde{\epsilon}_{\mathrm{eff}}^N& i\omega\tilde{\tilde{\mu}}_{\mathrm{eff}}^N \\
-i\omega\tilde{\tilde{\epsilon}}_{\mathrm{eff}}^N & \tilde{\mu}_{\mathrm{eff}}^N
\end{pmatrix}:= \Lambda  N^{-3a}c_n K_0 M_n(\epsilon_c,B).  
$$
Then, from \eqref{EN_decomp} and Proposition \ref{prop_pointsource},
we can see that
\begin{equation} \label{eq*}
\begin{pmatrix}
E^N(x)
\\
H^N(x)
\end{pmatrix}
= \begin{pmatrix}
E^{in}(x)
\\
H^{in}(x)
\end{pmatrix} + 
\frac{1}{N^3}\sum\limits_{j=1}^{N^3}  \omega G(x-z_j^N)
\begin{pmatrix}
\tilde{\epsilon}_{\mathrm{eff}}^N& i\omega\tilde{\tilde{\mu}}_{\mathrm{eff}}^N \\
-i\omega\tilde{\tilde{\epsilon}}_{\mathrm{eff}}^N & \tilde{\mu}_{\mathrm{eff}}^N
\end{pmatrix}
\begin{pmatrix}
E^{N}(z_j^N) \\ H^{N}(z_j^N)
\end{pmatrix}+O(\frac{ N^3\delta^4}{|\det A_n|})
\end{equation}
for $x\in\Omega\setminus (\cup_{j=1}^{N^3}z_j^N)$.
Assuming the homogenized limit $(E^{h},H^{h}):=\lim_{N\rightarrow \infty}({E^N},H^N)$ exists in some sense, we can easily expect that $(E^h,H^h)$ satisfies
\begin{equation}\label{hom_Lip}
\begin{pmatrix}
E^{h}
\\
H^{h}
\end{pmatrix}
= \begin{pmatrix}
E^{in}
\\
H^{in}
\end{pmatrix} + 
\int_{{\Omega}} \omega G(\cdot-z)  \begin{pmatrix}
\tilde{\epsilon}_{\mathrm{eff}}(z) & i\omega\tilde{\tilde{\mu}}_{\mathrm{eff}}(z) \\
-i\omega\tilde{\tilde{\epsilon}}_{\mathrm{eff}}(z) & \tilde{\mu}_{\mathrm{eff}}(z)
\end{pmatrix}
\begin{pmatrix}
E^{h}(z) \\ H^{h}(z) 
\end{pmatrix} dz \quad \mbox{in } {\Omega},
\end{equation}
where
\begin{align*} 
\begin{pmatrix}
\tilde{\epsilon}_{\mathrm{eff}}& i\omega\tilde{\tilde{\mu}}_{\mathrm{eff}} \\
-i\omega\tilde{\tilde{\epsilon}}_{\mathrm{eff}} & \tilde{\mu}_{\mathrm{eff}}
\end{pmatrix}    
 :=  \widetilde{V} \lim_{N\rightarrow \infty } \begin{pmatrix}
\tilde{\epsilon}_{\mathrm{eff}}^N& i\omega\tilde{\tilde{\mu}}_{\mathrm{eff}}^N \\
-i\omega\tilde{\tilde{\epsilon}}_{\mathrm{eff}}^N & \tilde{\mu}_{\mathrm{eff}}^N
\end{pmatrix}  .
\end{align*}
Straightforward but tedious computations show that the following compatibility condition holds:
\begin{equation} \label{compatible_1}
\ds
\frac{\ds \tilde{\tilde{\epsilon}}_{\mathrm{eff}} + \frac{\epsilon_m\beta_m}{1-\omega^2\epsilon_m\mu_m\beta_m^2}}
{\ds \tilde{\tilde{\mu}}_{\mathrm{eff}} + \frac{\mu_m\beta_m}{1-\omega^2\epsilon_m\mu_m\beta_m^2} } = \frac{\epsilon_m}{\mu_m}. 
\end{equation}

We can conclude by comparing \eqref{hom_Lip} with the Lippmann-Schwinger equation \eqref{eq:chiralgreen}, that the effective electric permittivity $\epsilon_{\mathrm{eff}}$ and magnetic permeability $\mu_{\mathrm{eff}}$ are determined by solving the 
system of three equations
\begin{align}
\tilde{\epsilon}_{\mathrm{eff}}& = \frac{\epsilon_{\mathrm{eff}}}{\epsilon_m(1-\omega^2\epsilon_{\mathrm{eff}}\mu_{\mathrm{eff}}\beta^2_{\mathrm{eff}})} - \frac{1}{1-\omega^2\epsilon_m\mu_m\beta_m^2},\nonumber \\
\tilde{\mu}_{\mathrm{eff}} &= \frac{\mu_{\mathrm{eff}}}{\mu_m (1-\omega^2\epsilon_{\mathrm{eff}}\mu_{\mathrm{eff}} \beta^2_{\mathrm{eff}})} - \frac{1}{1-\omega^2\epsilon_m\mu_m\beta_m^2},\nonumber \\
\tilde{\tilde{\epsilon}}_{\mathrm{eff}} &=  \frac{\epsilon_{\mathrm{eff}} \mu_{\mathrm{eff}}\beta_{\mathrm{eff}}}{\mu_m(1-\omega^2\epsilon_{\mathrm{eff}} \mu_{\mathrm{eff}}\beta^2_{\mathrm{eff}})}-\frac{\epsilon_m\beta_m}{1-\omega^2\epsilon_m\mu_m\beta_m^2}.
\label{tilde_nonlinear_eqns}
\end{align}
Here, $\beta_{\mathrm{eff}}$ is the effective chiral admittance.   Note that, due to the compatibility condition \eqref{compatible_1}, the following relation immediately holds:
$$
\tilde{\tilde{\mu}}_{\mathrm{eff}} =  \frac{\epsilon_{\mathrm{eff}}\mu_{\mathrm{eff}} \beta_{\mathrm{eff}}}{\epsilon_m(1-\omega^2\epsilon_{\mathrm{eff}}\mu_{\mathrm{eff}}\beta^2_{\mathrm{eff}})} - \frac{\mu_m\beta_m}{1-\omega^2\epsilon_m\mu_m\beta_m^2}.
$$
In fact, the equation \eqref{tilde_nonlinear_eqns} is uniquely solvable.
One can easily check that
\begin{align}
\epsilon_{\mathrm{eff}} &= \epsilon_m \left( (\tilde{\epsilon}_{\mathrm{eff}} + \tilde{\gamma}_m) -\omega^2 (\tilde{\tilde{\epsilon}}_{\mathrm{eff}} + \epsilon_m \beta_m \tilde{\gamma}_m)\frac{\tilde{\tilde{\mu}}_{\mathrm{eff}} + \mu_m\beta_m\tilde{\gamma}_m}{\tilde{\mu}_{\mathrm{eff}} + \tilde{\gamma}_m}\right), 
\nonumber
\\
\mu_{\mathrm{eff}} &= \mu_m \left( (\tilde{\mu}_{\mathrm{eff}} + \tilde{\gamma}_m) -\omega^2 (\tilde{\tilde{\mu}}_{\mathrm{eff}} + \mu_m \beta_m \tilde{\gamma}_m)\frac{\tilde{\tilde{\epsilon}}_{\mathrm{eff}} + \epsilon_m\beta_m\tilde{\gamma}_m}{\tilde{\epsilon}_{\mathrm{eff}} + \tilde{\gamma}_m}\right), 
\label{eps_mu_eff_formula}
\end{align}
where
$$
\tilde{\gamma}_m = \frac{1}{1-k_m^2\beta_m^2}.
$$
Therefore, $(E^h,H^h)$ satisfies 
\begin{equation}\label{eq:chiralmaxwell_hom}
\left\{
\begin{aligned}
&\nabla\times E^h = \gamma_h^2\beta_h E^h+i\omega\mu_h \left(\frac{\gamma_h}{k_h}\right)^2 H^h \quad\mathrm{in}~\Omega,\\
&\nabla\times H^h = \gamma_h^2\beta_h H^h-i\omega\epsilon_h\left(\frac{\gamma_h}{k_h}\right)^2 E^h \quad\mathrm{in}~\Omega,
\end{aligned}
\right.
\end{equation}
with the homogenized material parameters
$$
\epsilon_h = \begin{cases}
 \epsilon_{\mathrm{eff}} , &\quad \mbox{in }\Omega
 \\
 \epsilon_m, &\quad \mbox{in }\mathbb{R}^3\setminus\Omega
 \end{cases},
 \quad
\mu_h = \begin{cases}
 \mu_{\mathrm{eff}} , &\quad \mbox{in }\Omega
 \\
 \mu_m, &\quad \mbox{in }\mathbb{R}^3\setminus\Omega
 \end{cases}. 
$$
The other parameters $\beta_h, \gamma_h, $ and $k_h$ are defined similarly.

\subsection{Double-negative effective properties}

Here, we show that  the effective properties $\epsilon_{\mathrm{eff}}$ and $\mu_{\mathrm{eff}}$ can be both negative.

Let us first consider the resonant behavior of the matrix
$$
\begin{pmatrix}
\tilde{\epsilon}_{\mathrm{eff}}^N& i\omega\tilde{\tilde{\mu}}_{\mathrm{eff}}^N \\
-i\omega\tilde{\tilde{\epsilon}}_{\mathrm{eff}}^N & \tilde{\mu}_{\mathrm{eff}}^N
\end{pmatrix}, 
$$
when $\epsilon_c$ is close to $\epsilon_{c,n}^*$.
Straightforward  but tedious computations show that each of components has the following behavior with respect to $\epsilon_c$:
\begin{align}
\tilde{\epsilon}_{\mathrm{eff}}^N &= \frac{\Lambda N^{-3\alpha} c_n \epsilon_m }{(1/2-\lambda_n)^3(1-\beta_m^2 k_m^2(1/2-\lambda_n))^2 } \Big(\frac{-1}{\epsilon_c - \epsilon_{c,n}^*} + O(1)\Big),\nonumber
\\
\tilde{\mu}_{\mathrm{eff}}^N &= \frac{ \Lambda N^{-3\alpha}c_n\epsilon_m k_m^2\beta_m^2 }{(1/2-\lambda_n)(1-\beta_m^2 k_m^2(1/2-\lambda_n))^2 } \Big(\frac{-1}{\epsilon_c - \epsilon_{c,n}^*} + O(1)\Big),\nonumber
\\
\tilde{\tilde{\epsilon}}_{\mathrm{eff}}^N &= \frac{ \Lambda N^{-3\alpha}c_n\epsilon_m^2 \beta_m }{(1/2-\lambda_n)^2(1-\beta_m^2 k_m^2(1/2-\lambda_n))^2 } \Big(\frac{-1}{\epsilon_c - \epsilon_{c,n}^*} + O(1)\Big),\nonumber
\\
\tilde{\tilde{\mu}}_{\mathrm{eff}}^N &= \frac{ \Lambda N^{-3\alpha}c_n\epsilon_m \mu_m \beta_m }{(1/2-\lambda_n)^2(1-\beta_m^2 k_m^2(1/2-\lambda_n))^2 } \Big(\frac{-1}{\epsilon_c - \epsilon_{c,n}^*} + O(1)\Big).
\label{asymp_ep_mu_tilde_tilde2}
\end{align}
Here $O(1)$ means that the remainder  does not diverge for any $\epsilon_c$. 

Now we turn to the effective properties $\epsilon_{\mathrm{eff}}(y)$ and $\mu_{\mathrm{eff}}(y)$ for $y\in\Omega$. For the sake of simplicity of presentation, we assume that at $y$, $\widetilde{V}(y)=1$. 
By applying the asymptotics \eqref{asymp_ep_mu_tilde_tilde2} to \eqref{eps_mu_eff_formula}, one can check that each of $\epsilon_{\mathrm{eff}}$ and $\mu_{\mathrm{eff}}$ has a removable singularity at $\epsilon_c = \epsilon_{c,n}^*$. In fact, $\epsilon_{\mathrm{eff}}$ (or $\mu_{\mathrm{eff}}$) diverges only when $\epsilon_c$ satisfies $\tilde{\mu}_{\mathrm{eff}} + \tilde{\gamma}_m=0$ (respectively  $\tilde{\epsilon}_{\mathrm{eff}} + \tilde{\gamma}_m=0$).
Let $\epsilon_{c,n}^*[\epsilon_{\mathrm{eff}}]$ (and $\epsilon_{c,n}^*[\mu_{\mathrm{eff}}]$) be the value of $\epsilon_c$ at which $\epsilon_{\mathrm{eff}}$ (respectively $\mu_{\mathrm{eff}}$) diverges.
Then, it can be easily checked that, for large $N$,
\begin{align}
\epsilon_{c,n}^*[\epsilon_{\mathrm{eff}}] \approx \epsilon_{c,n}^* + 
\frac{ \epsilon_m k_m^2\beta_m^2 }{(1/2-\lambda_n)(1-\beta_m^2 k_m^2(1/2-\lambda_n))^2 } \tilde{\gamma}_m^{-1}c_n\Lambda N^{-3\alpha} 
\\
\epsilon_{c,n}^*[\mu_{\mathrm{eff}}] \approx \epsilon_{c,n}^* + 
\frac{ \epsilon_m  }{(1/2-\lambda_n)^3(1-\beta_m^2 k_m^2(1/2-\lambda_n))^2 } \tilde{\gamma}_m^{-1}c_n\Lambda N^{-3\alpha} 
\end{align}
Note that $$\epsilon_{c,n}^*[\mu_{\mathrm{eff}}] > \epsilon_{c,n}^*[\epsilon_{\mathrm{eff}}].$$

Now we choose $\epsilon_c$ such that $\epsilon_c$ is slightly above the two resonant permittivities. More precisely, we assume the following.

\begin{asump}\label{asump_epc2}
Let $0<s<1$ be given and assume that
\begin{equation}\label{epc_choice}
\epsilon_c = \epsilon_{c,n}^* + s^{-1}
\frac{ \Lambda N^{-3\alpha}c_n \epsilon_m  }{(1/2-\lambda_n)^3(1-\beta_m^2 k_m^2(1/2-\lambda_n))^2 } \tilde{\gamma}_m^{-1}. 
\end{equation}
If $s$ is very close to one, $\epsilon_c $ is slightly above the resonant permittivity $\epsilon_{c,n}^*[\mu_{\mathrm{eff}}]$.
\end{asump}

By substituting \eqref{epc_choice} into \eqref{asymp_ep_mu_tilde_tilde2} and then taking limit $N\rightarrow \infty$, we have 
\begin{align*}
\tilde{\epsilon}_{\mathrm{eff}}& = -s\tilde{\gamma}_m,
\\
\tilde{\mu}_{\mathrm{eff}} &= -s\tilde{\gamma}_m (k_m\beta_m)^2(1/2-\lambda_n)^2
, \\
\tilde{\tilde{\epsilon}}_{\mathrm{eff}} &= -s\tilde{\gamma}_m\epsilon_m \beta_m (1/2-\lambda_n),
\\
\tilde{\tilde{\mu}}_{\mathrm{eff}} &= -s\tilde{\gamma}_m\mu_m \beta_m (1/2-\lambda_n).
\end{align*}
So, using \eqref{eps_mu_eff_formula}, we finally get the formulas for the effective material parameters as follows:
\begin{align}
\epsilon_{\mathrm{eff}} &= \epsilon_m \tilde{\gamma}_m \Big( (1-s) -\omega^2 \big(1-s(1/2-\lambda_n)\big)\frac{\mu_m\beta_m(1-s(1/2-\lambda_n))}{1- s (k_m\beta_m)^2 (1/2-\lambda_n)^2}\Big),
\nonumber
\\
\mu_{\mathrm{eff}} &= \mu_m \tilde{\gamma}_m \Big( 1-s(k_m\beta_m)^2(1/2-\lambda_n)^2
-  \frac{1 }{1-s} 
 (k_m\beta_m)^2 \Big( 1-s(1/2-\lambda_n))^2 
\Big).
\label{eps_mu_eff_formula_final}
\end{align}

Now we can prove that the above effective parameters are both negative as shown in the following theorem.

\begin{theorem} {\rm(Double-negative property)}
Suppose that the permittivity $\epsilon_c$ of the plasmonic particle is given as in Assumption \ref{asump_epc2}.
Then, the effective parameters $\epsilon_{\mathrm{eff}}$ and $\mu_{\mathrm{eff}}$ of the homogenized equation are both negative provided $s$ is sufficiently close to one.

\end{theorem}
\proof
Since $0<s<1$, $|\lambda_n|<1/2$ and $k_m\beta_m <1$, the conclusion immediately follows from \eqref{eps_mu_eff_formula_final}. 
\qed

\bigskip
\begin{rmk}\label{rmk:random}
Although we assume the shape of the particle is symmetric, our result can be extended to the case of arbitrary shaped particles. Suppose that the particles are randomly oriented. Then the average of $\langle\nu,\phi_n\rangle_*\langle\nu,\phi_n\rangle_*^\top$ over the orientation of the particle becomes a diagonal matrix as in the symmetric case \eqref{tensor_cn}. Therefore, we obtain in exactly the same manner 
negative effective permittivity and negative effective permeability for frequencies near the resonant permittivity $\epsilon_{c,n}^*$.
\end{rmk}

\begin{figure*}
\begin{center}
\epsfig{figure=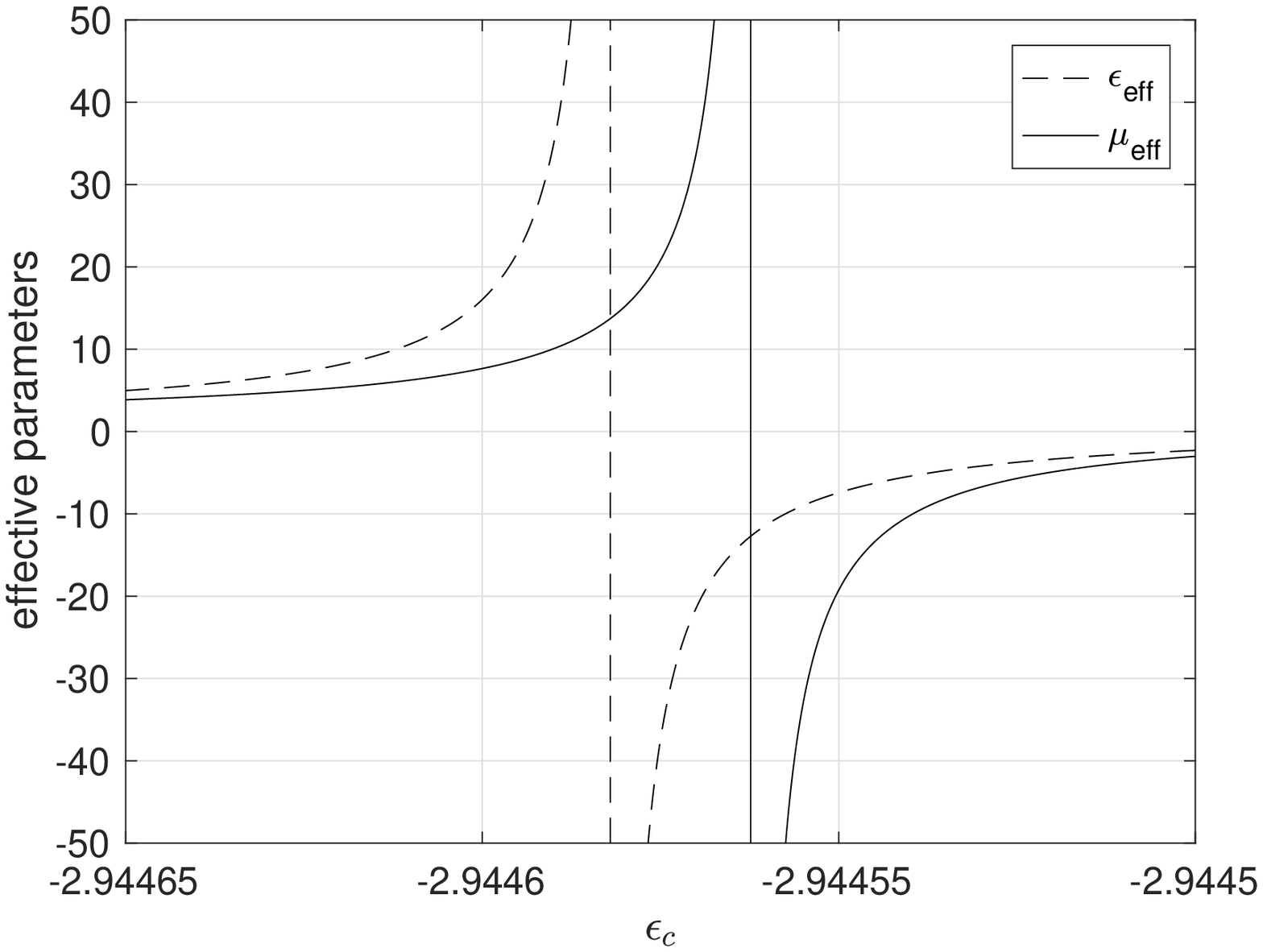,width=7cm}
\epsfig{figure=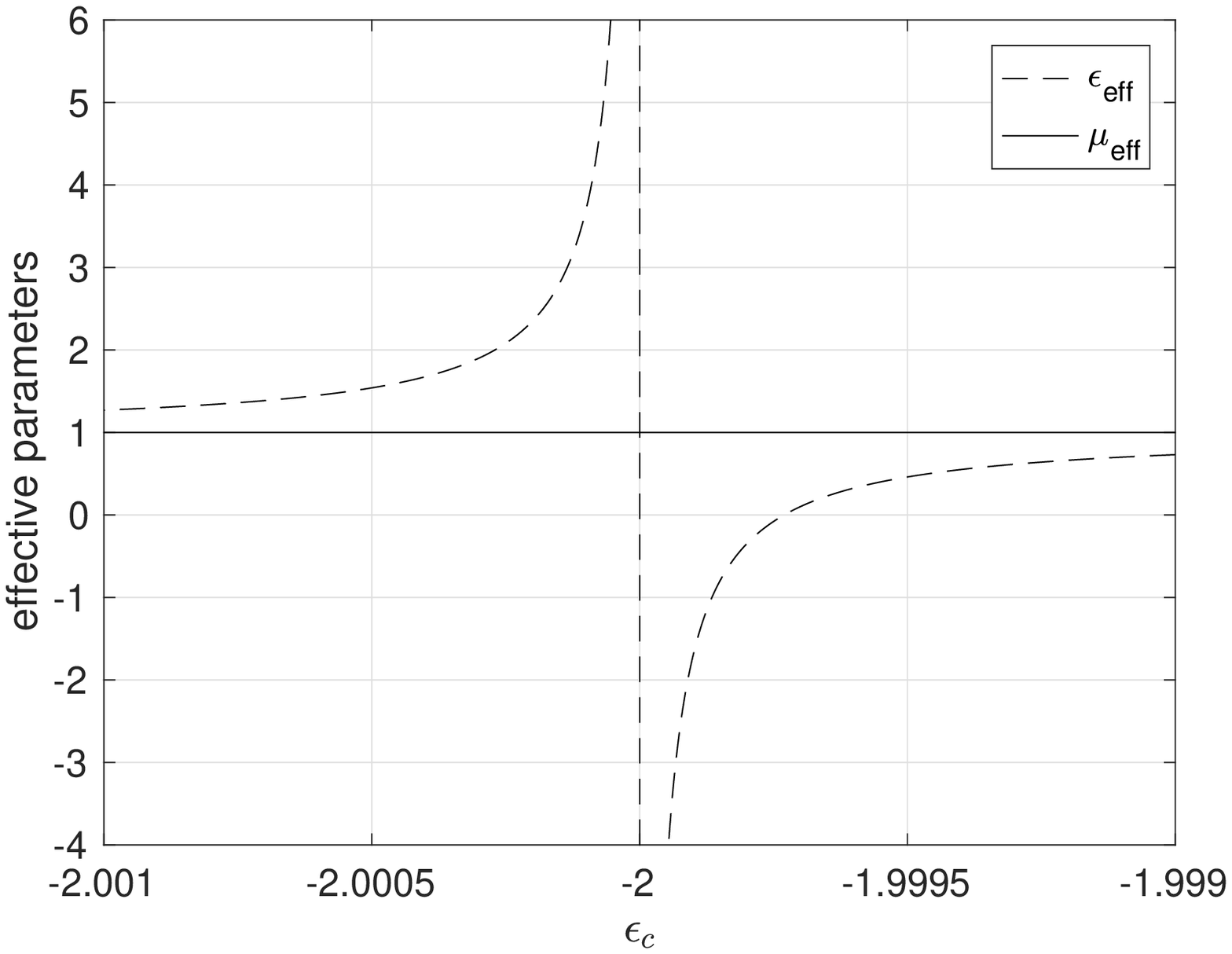,width=7cm}
\end{center}
\caption{The effective properties of the homogenized media.  A chiral media $\beta_m\neq 0$ (left),   a non-chiral media $\beta_m= 0$ (right).}
\label{fig:eff}
\end{figure*}

\begin{rmk}
We provide a numerical example in Figure \ref{fig:eff}. We plot the effective parameters $\epsilon_{\mathrm{eff}}$ and $\mu_{\mathrm{eff}}$ as functions of $\epsilon_c$. We set $\omega=1, \epsilon_m=1, \mu_m=1$ $\Lambda=3, N=125$ and $a=0.965$. In the left figure, we use $\beta_m = 1.09$. Clearly, both the effective parameters are negative near $\epsilon_c = -2.94455$. In the right figure, we change $\beta_m$ as $\beta_m=0$, which means that there is no chirality. In this case, only $\epsilon_{\mathrm{eff}}$ is resonant but  $\mu_{\mathrm{eff}}$ remains as one. This shows the importance of the chirality to achieve the double-negative metamaterial.
\end{rmk}

\begin{rmk}
The resonance frequency can be determined by the Drude model
$$
\epsilon_c(\omega)= 1 - \frac{\omega_p^2}{\omega^2 + i \tau \omega },
$$
where $\omega_p$ and $\tau$ are two given positive constants. 
\end{rmk}

\begin{rmk}
For simplicity, we assume that $\epsilon_c$ is real. The analysis in this section applies to the case where $\mathrm{Im} \epsilon_c$ is sufficiently small. 
\end{rmk}

\begin{rmk} 
By using the same approach as in \cite{haisima,figari,papanicolaou},  we provide under some assumptions on the distribution of the plasmonic particles a  justification of  the derivation of  the effective medium parameters. See Appendix \ref{appenC}. 
\end{rmk}

\section{Concluding remarks} 

In this paper, we have first derived an asymptotic expansion of the scattered electromagnetic fields by a small plasmonic dielectric nanoparticle in a chiral medium. We have shown that the plasmonic particle can be approximated by the sum of a resonant electric dipole and a resonant magnetic dipole. We have also characterized these resonant frequencies in terms of the chirality admittance of the background medium and the material parameters and the shape of the particle. Then we have obtained an effective medium theory for materials consisting of a large number of plasmonic nanoparticles embedded in a chiral background medium. We have shown that the dielectric plasmonic particles contribute to both the effective electric permittivity and the effective magnetic permeability. Finally, we have proved that both the effective electric permittivity and magnetic permeability can be negative near some resonant frequencies.

\appendix

\section{Explicit computation for a ball}\label{app_ball}

Suppose that $B$ is the unit ball. In this case, we are able to write out its polarization tensor $M(\epsilon_c,B)$ defined in (\ref{eq:polarizationdefine}) explicitly.
In order to do so, we compute the tensor $\sum\limits_{l=-1}^1\langle\nu,Y_1^l\rangle_*\langle\nu,Y_1^l\rangle_*^\top$.
 The following lemma from \cite{HA3} will be required. 
\begin{lemma} \label{lemball}
	For $n=0,1,\cdots$, we have
	$$
	\mathcal{K}_B^*[Y_n^l]=\frac{1}{2(2n+1)}Y_n^l(\hat{x}),\quad |x|=1,l=-n,\cdots,n,
	$$
	where $\hat{x}=x/|x|$ and $(Y_n^l)_{l=-n,\cdots,n}$ are the orthonormal spherical harmonics of degree $n$ and order $l$. Moreover, 
	\begin{align*}
\langle\nu, Y_n^l\rangle_* &= -\int_{\partial\Omega} \mathcal{S}_B [Y_n^l]\nu  \mathrm{d}\sigma(x) \\
&= -\int_{\partial B} -\frac{1}{2n+1} Y_n^l(\hat{x})\hat{x}\mathrm{d}\sigma(x) \\
&= \frac{1}{2n+1}\int_{\partial B} Y_n^l(\hat{x})\hat{x}\mathrm{d}\sigma(\hat{x}).
\end{align*}
	
\end{lemma} 
 Since $\hat{x}=(\sin\theta\cos\varphi, \sin\theta\sin\varphi,\cos\theta)$, and by definition of the spherical harmonic functions,
$$
Y_n^l(\theta,\varphi):=\sqrt{\frac{2n+1}{4\pi}\frac{(n-|l|)!}{(n+|l|)!}}P_n^{|l|}(\cos\theta)e^{il\varphi},
$$
with $P_n^{|l|}$ being the associated Legendre polynomial of degree $n$ and order $|l|$, 
we have
\begin{align*}
Y_1^{-1}(\theta,\varphi) &= \frac{1}{2}\sqrt\frac{3}{2\pi}P_1^1(\cos\theta)e^{-i\varphi} = \frac{1}{2}\sqrt\frac{3}{2\pi}\sin\theta e^{-i\varphi}, \\
Y_1^1(\theta,\varphi) &= \frac{1}{2}\sqrt\frac{3}{2\pi}P_1^{1}(\cos\theta)e^{i\varphi} = \frac{1}{2}\sqrt\frac{3}{2\pi}\sin\theta e^{i\varphi}, \\
Y_1^0(\theta,\varphi) &= \frac{1}{2}\sqrt\frac{3}{\pi}P_1^0(\cos\theta)=\frac{1}{2}\sqrt\frac{3}{\pi}\cos\theta.
\end{align*}
Consequently,
\begin{align*}
\sin\theta\cos\varphi&=\sqrt\frac{2\pi}{3}(\bar{Y}_1^{-1} + \bar{Y}_1^1), \\
\sin\theta\sin\varphi&=i\sqrt\frac{2\pi}{3}(\bar{Y}_1^1-\bar{Y}_1^{-1}), \\
\cos\theta &= 2\sqrt\frac{\pi}{3}\bar{Y_1^0},
\end{align*}
where $\bar{Y}_n^l$ is  the complex conjugate of $Y_n^l$. Since $\{Y_n^l\}$ is an orthogonal basis of $L^2(\partial B)$, the infinite sum in \eqref{eq:polarizationdefine} is actually finite, and among all the terms only $\langle \nu,Y_1^{-1}\rangle_*$, 
$\langle \nu,Y_1^0\rangle_*$,  and $\langle \nu,Y_1^1\rangle_*$ are nonzero. We can calculate that
\begin{align*}
\sum\limits_{l,n}\langle\nu,Y_n^l\rangle_*\langle\nu,Y_n^l\rangle_*^\top
&=\sum\limits_{l=-1}^1\langle\nu,Y_1^l\rangle_*\langle\nu,Y_1^l\rangle_*^\top
\\
&=\frac{1}{9}\left(\begin{pmatrix}
-\sqrt\frac{2\pi}{3} \\ i\sqrt\frac{2\pi}{3} \\ 0\end{pmatrix}\begin{pmatrix}
-\sqrt\frac{2\pi}{3} & -i\sqrt\frac{2\pi}{3} & 0
\end{pmatrix} + \begin{pmatrix}
0 \\ 0 \\ 2\sqrt\frac{\pi}{3}
\end{pmatrix}\begin{pmatrix} 0 & 0 & 2\sqrt\frac{\pi}{3}\end{pmatrix} \right) \\
&\quad + \frac{1}{9} \left( \begin{pmatrix}
-\sqrt\frac{2\pi}{3} \\ -i\sqrt\frac{2\pi}{3} \\ 0\end{pmatrix}\begin{pmatrix}
-\sqrt\frac{2\pi}{3} & i\sqrt\frac{2\pi}{3} & 0
\end{pmatrix}\right)\\
&=\frac{4\pi}{27}I.
\end{align*}

\section{Layer potentials for electromagnetic waves in a chiral medium}\label{subsec:layer_chiral}

In this appendix, we briefly review the results in \cite{ima} concerning the layer potential techniques for electromagnetic scattering by the particle $\Omega$ in a chiral medium.

For $s=\pm1/2$, let $H^s(\p \Om)$ denote the usual Sobolev space of order $s$ on $\p \Om$ and let $$H^{s}_T(\p \Om) = \left\{\varphi \in \big(H^{s}(\p \Om)\big)^3, \nu\cdot\varphi=0  \right\}.$$ Let $\nabla_{\p \Om}$, $\nabla_{\p \Om}\cdot$ and $\Delta_{\p \Om}$ denote the surface gradient, surface divergence and Laplace-Beltrami operator respectively and define the vectorial and scalar surface curl by $\vec{\text{curl}}_{\p \Om}\varphi = -\nu\times\nabla_{\p \Om}\varphi$ for $\varphi \in H^{\f{1}{2}}(\p \Om)$ and $\text{curl}_{\p \Om}\varphi = -\nu \cdot(\nabla_{\p \Om} \times \varphi)$ for $\varphi \in H^{-\f{1}{2}}_T(\p \Om)$, respectively.
We introduce the following functional space:
\beas
H^{-\f{1}{2}}_T(\text{div},\p \Om) &=& \left\{ \varphi \in H^{-\f{1}{2}}_T(\p \Om), \nabla_{\p \Om}\cdot\varphi\in H^{-\f{1}{2}}(\p \Om) \right\}.
\eeas

We introduce the boundary layer potentials by
$$\begin{array}{l}
\ds \vec{\mathcal{S}}_{\Om}^{k} [\varphi](x) = \int_{\p \Om} g^k(x-y) \varphi(y) d\sigma(y),  \quad x 
\in \R^3,\\
\nm \ds
\mathcal{S}_{\Om}^{k} [\psi](x) = \int_{\p \Om} g^k(x-y) \psi(y) d\sigma(y) \quad \mbox{for a scalar function } \psi \in H^{-\f{1}{2}}(\p \Om) \mbox{ and } x \in \R^3,
\\
\nm \ds
\mathcal{M}_{\Om}^{k} [\varphi](x) = \int_{\p \Om} \nu(x)\times\nabla_x\times g^k(x-y)\varphi(y) d\sigma(y),  \quad x \in \p \Om,
\\
\nm \ds
\mathcal{L}_{\Om}^{k} [\varphi](x) = \nu(x)\times \bigg(k^2\vec{\mathcal{S}}_{\Om}^{k}[\varphi](x) + \nabla\mathcal{S}_\Om^k[\nabla_{\p \Om}\cdot\varphi](x)\bigg) ,   \quad x \in \p \Om.
\end{array}$$
 We also introduce the notations
\begin{align*}
&\ds\alpha_{1,c} = \frac{1}{i\tau_c}, \quad \alpha_{2,c} = -i \tau_c,
\quad \ds \alpha_{1,m} = \frac{1}{i\tau_m}, \quad \alpha_{2,m} = -i \tau_m, 
\end{align*}
where $
\tau_c = \sqrt{\mu_m/\epsilon_c}$ and $\tau_m = \sqrt{\mu_m/\epsilon_m}
$.

We consider the Bohren decomposition of $(E,H)$ into Beltrami fields,i.e., 
\begin{equation}\label{bohren_int}
E = Q_1 +\alpha_{2,c}  Q_2, \quad H = \alpha_{1,c}Q_1 + Q_2 \quad \mbox{in }\Omega.
\end{equation}
Similarly,
\begin{equation}\label{bohren_ext}
E^{sc} = Q_1 +\alpha_{2,m} Q_2, \quad H^{sc} = \alpha_{1,m}Q_1 + Q_2 \quad \mbox{in } \R^3\setminus\overline{\Omega}.
\end{equation}
We can see that they satisfy the vector Helmholtz equations as
$$
\begin{cases}
(\Delta + \gamma_{j,c}^2) Q_j = 0,&\quad \mbox{in }\Omega,
\\
(\Delta + \gamma_{j,m}^2) Q_j = 0,&\quad  \mbox{in }\R^3\setminus\overline{\Omega},
\end{cases}
$$
where
$$
\gamma_{j,m} = \frac{k_m}{1 + (-1)^j k_m \beta_m}, \quad \gamma_{j,c} = \omega \sqrt{\epsilon_c \mu_m}, \quad j=1,2.
$$


We define the operator $\mathcal{Q}_j$,  for $(\varphi_1,\varphi_2)\in H^{-\f{1}{2}}_T(\text{div},\p \Om)$, by
\begin{align*}
\mathcal{Q}_{\Omega,j} 
\begin{bmatrix} \varphi_1 \\ \varphi_2 \end{bmatrix}
 = 
\begin{cases}
\big((-1)^{j+1}\gamma_{j,c} \nabla \times \vec{\mathcal{S}}_{\Om}^{\gamma_{j,c}} + 
\nabla \times \nabla \times \vec{\mathcal{S}}_{\Om}^{\gamma_{j,c}}\big)[\varphi_j]  &\quad \mbox{in }\Omega, 
\\[0.5em]
\big((-1)^{j+1}\gamma_{j,c} \nabla \times  \vec{\mathcal{S}}_{\Om}^{\gamma_{j,m}} + 
\nabla \times \nabla \times  \vec{\mathcal{S}}_{\Om}^{\gamma_{j,m}}\big) [\zeta_{j1}\varphi_1 + \zeta_{j2}\varphi_2] &\quad \mbox{in }\R^3\setminus\overline{\Omega}
\end{cases}
\end{align*}
with $\zeta_{ij},i,j=1,2$, given by
\begin{align}
\zeta_{11} &= \frac{1}{2}(1+\frac{\tau_m}{\tau_c}),
\quad
\zeta_{12} = \frac{i}{2}({\tau_c}-{\tau_m}),
\\
\zeta_{21} &= \frac{i}{2}(\f{1}{\tau_m}-\f{1}{\tau_c}),
\quad
\zeta_{22} = \frac{i}{2}(1+\f{\tau_c}{\tau_m}).
\end{align}
Then the solution $Q_j,$ for $j=1,2$, can be represented as
\begin{equation}\label{Qj_rep}
Q_j = \mathcal{Q}_{\Omega,j} \begin{bmatrix} \psi_1 \\ \psi_2 \end{bmatrix},
\end{equation}
where  $(\psi_1,\psi_2)$ is the solution of the integral equation
\begin{equation}\label{eq**}
(\mathcal{J}_{\Omega} + \mathcal{C}_\Omega) \begin{pmatrix}
 \psi_1 \\ \psi_2
 \end{pmatrix}
 = \begin{pmatrix}
 \nu\times E^{in}|_{\p \Om} \\ \nu\times H^{in}|_{\p \Om}
 \end{pmatrix}.
\end{equation}
Here, the operator $\mathcal{J}_\Omega$ is given by
\begin{align*}
(\mathcal{J}_\Omega)_{11} &= -\zeta_{11}\mathcal{L}^{\gamma_{1,m}}_\Omega
+ \mathcal{L}^{\gamma_{1,c}}_\Omega - \zeta_{21} \alpha_{2,m}
\mathcal{L}^{\gamma_{2,m}}_\Omega
-\frac{1}{2} (\zeta_{11} \gamma_{1,m} + \gamma_{1,c} -\zeta_{21}\alpha_{2,m}\gamma_{2,m} )I,
\\
(\mathcal{J}_\Omega)_{12} &= -\zeta_{11}\mathcal{L}^{\gamma_{1,m}}_\Omega
+ \alpha_{2,m}\mathcal{L}^{\gamma_{2,m}}_\Omega - \zeta_{22} \alpha_{2,m}
\mathcal{L}^{\gamma_{2,m}}_\Omega
-\frac{1}{2} (\zeta_{12} \gamma_{1,m} -\alpha_{2,c} \gamma_{2,c} -\zeta_{22}\alpha_{2,m}\gamma_{2,m} )I,
\\
(\mathcal{J}_\Omega)_{21} &= -\zeta_{21}\mathcal{L}^{\gamma_{2,m}}_\Omega
+ \alpha_{1,m}\mathcal{L}^{\gamma_{1,m}}_\Omega - \zeta_{11} \alpha_{1,m}
\mathcal{L}^{\gamma_{1,m}}_\Omega
-\frac{1}{2} (-\zeta_{21} \gamma_{2,m} + \alpha_{1,c}\gamma_{1,c} +\zeta_{11}\alpha_{1,m}\gamma_{1,m} )I,
\\
(\mathcal{J}_\Omega)_{22} &= -\zeta_{22}\mathcal{L}^{\gamma_{2,m}}_\Omega
+ \mathcal{L}^{\gamma_{2,c}}_\Omega
 - \zeta_{12} \alpha_{1,m}\mathcal{L}^{\gamma_{1,m}}_\Omega
-\frac{1}{2} (\zeta_{22} \gamma_{2,m} - \gamma_{2,c} +\zeta_{12}\alpha_{1,m}\gamma_{1,m} )I,
\end{align*}
and the operator $\mathcal{C}_\Omega$ is given by
\begin{align*}
(\mathcal{C}_\Omega)_{11} &= 
 -\zeta_{11} \gamma_{1,m}\mathcal{M}_\Omega^{\gamma_{1,m}} - \gamma_{1,c}\mathcal{M}_\Omega^{\gamma_{1,c}} +\zeta_{21}\alpha_{2,m}\gamma_{2,m}\mathcal{M}_\Omega^{\gamma_{2,m}},
\\[0.5em]
(\mathcal{C}_\Omega)_{12} &= 
 -\zeta_{12} \gamma_{1,m}\mathcal{M}_\Omega^{\gamma_{1,m}} +\alpha_{2,c} \gamma_{2,c} \mathcal{M}_\Omega^{\gamma_{2,c}} +\zeta_{22}\alpha_{2,m}\gamma_{2,m} \mathcal{M}_\Omega^{\gamma_{2,m}},
\\[0.5em]
(\mathcal{C}_\Omega)_{21} &= 
 \zeta_{21} \gamma_{2,m}\mathcal{M}_\Omega^{\gamma_{2,m}} - \alpha_{1,c}\gamma_{1,c}\mathcal{M}_\Omega^{\gamma_{1,c}} -\zeta_{11}\alpha_{1,m}\gamma_{1,m}\mathcal{M}_\Omega^{\gamma_{1,m}},
\\[0.5em]
(\mathcal{C}_\Omega)_{22} &= 
 -\zeta_{22} \gamma_{2,m}\mathcal{M}_\Omega^{\gamma_{2,m}} + \gamma_{2,c}\mathcal{M}_\Omega^{\gamma_{2,c}} -\zeta_{12}\alpha_{1,m}\gamma_{1,m}\mathcal{M}_\Omega^{\gamma_{1,m}}.
\end{align*}
The operator $\mathcal{C}_\Omega: H^{-\f{1}{2}}_T(\text{div},\p \Om) \rightarrow H^{-\f{1}{2}}_T(\text{div},\p \Om)$ is compact.

Let
$$
\mathcal{Q}^E_\Omega = \mathcal{Q}_{\Omega,1} + \alpha_{2,c}  \mathcal{Q}_{\Omega,2}, 
\quad
\mathcal{Q}^H_\Omega =  \alpha_{1,c}\mathcal{Q}_{\Omega,1} + \mathcal{Q}_{\Omega,2}.
$$
In view of \eqref{bohren_int},  \eqref{bohren_ext}  and  \eqref{Qj_rep}, the operator $\mathcal{Q}^E_\Omega$ (or $\mathcal{Q}^H_\Omega$) maps density functions on $\p\Omega$ to the corresponding electric field (respectively,  the magnetic field).

 \section{Justification of the homogenization procedure} \label{appenC}
In this appendix,  we provide a justification of the point interaction approximation for deriving the effective medium parameters. We make the following assumptions.
\begin{asump} \label{asumpfinal}
\begin{itemize}
\item[(i)] The function $\widetilde{V}$ defined in Assumption \ref{assumptionV1} belongs to $\mathcal{C}^2_0(\Omega)$. 
 \item[(ii)]  (\ref{assumptionV2}) holds for functions $(f,g) \in X$, where $X:= \mathcal{C}^{0,\alpha}({\Omega})^3 \times \mathcal{C}^{0,\alpha}({\Omega})^3$ with $0< \alpha \leq 1$.
 \end{itemize} 
 \end{asump}
We define the operator $\mathcal{T}: X  \rightarrow X$ by
\begin{equation}\label{operatorT}
\mathcal{T}\begin{pmatrix}
u \\ v
\end{pmatrix}
=
\omega \int_{{\Omega}}G(\cdot-z) 
\begin{pmatrix}
\tilde{\epsilon}_{\mathrm{eff}}(z) & i\omega\tilde{\tilde{\mu}}_{\mathrm{eff}}(z) \\
-i\omega\tilde{\tilde{\epsilon}}_{\mathrm{eff}}(z) & \tilde{\mu}_{\mathrm{eff}}(z)
\end{pmatrix}
\begin{pmatrix}
u(z) \\ v(z)
\end{pmatrix}\; dz.
\end{equation}
Then, the Lippmann-Schwinger equation can be written as
$$
(I-\mathcal{T} ) \begin{pmatrix}
E^{h} \\ H^{h}
\end{pmatrix}
=\begin{pmatrix}
E^{in} \\ H^{in}
\end{pmatrix}.
$$

We assume that the homogenized problem \eqref{eq:chiralmaxwell_hom} is well-posed. More precisely, we make the following assumption.
\begin{asump}\label{asump-6}
For given material parameters $\epsilon_{\mathrm{eff}}$ and $\mu_{\mathrm{eff}}$ with negative real parts and
an incident field $(E^{in},H^{in})$, there exists a unique solution to \eqref{eq:chiralmaxwell_hom} such that $(E-E^{in},H-H^{in})$ satisfies the Silver-M\"uller radiation condition at infinity. 
\end{asump}

\begin{lemma}
The following statements are equivalent: 
\begin{enumerate}
\item[(i)]
There exists a unique solution to the differential equation \eqref{eq:chiralmaxwell_hom} such that $(E-E^{in},H-H^{in})$ satisfies the Silver-M\"uller radiation condition at infinity;

\item[(ii)]
There exists a unique solution $(E,H) \in X$ to the Lippmann-Schwinger equation
$$
(I-\mathcal{T} ) \begin{pmatrix}
E \\ H
\end{pmatrix}
=\begin{pmatrix}
E^{in} \\ H^{in}
\end{pmatrix}, 
$$
where $\mathcal{T}$ is given by \eqref{operatorT}.
\end{enumerate}

\end{lemma}

Note that, under Assumptions \ref{asumpfinal} (i) and \ref{asump-6}, the operator $\mathcal{T}$ is Fredholm of index zero on the set of functions $(u,v) \in X$ such that \cite{CK}
$$
\nabla \cdot \bigg( \bigg[   
\begin{pmatrix}
\tilde{\epsilon}_{\mathrm{eff}} (z) & i\omega\tilde{\tilde{\mu}}_{\mathrm{eff}}(z) \\
-i\omega\tilde{\tilde{\epsilon}}_{\mathrm{eff}} (z) & \tilde{\mu}_{\mathrm{eff}}(z)
\end{pmatrix} + \begin{pmatrix} 
 \frac{1}{1-\omega^2\epsilon_m\mu_m\beta_m^2}  &  i \omega \frac{\mu_m\beta_m}{1-\omega^2\epsilon_m\mu_m\beta_m^2}\\
- i \omega \frac{\epsilon_m\beta_m}{1-\omega^2\epsilon_m\mu_m\beta_m^2} &  \frac{1}{1-\omega^2\epsilon_m\mu_m\beta_m^2} \end{pmatrix}\bigg] 
\begin{pmatrix}
u(z) \\ v(z)
\end{pmatrix}   \bigg) = 0 \quad \mbox{for } z \in {\Omega}.
$$

Let us introduce a regularized operator $\mathcal{T}_\eta$ of $\mathcal{T}$ by replacing $g^k$ with 
$$
	g_{\eta}^k(x,y) = \frac{e^{ik|x-y|}}{4\pi|x-y|+\eta},\qquad \eta>0.
$$
It is clear that the operator $\mathcal{T}_\eta$ is compact in $X$. Using Fredholm's theory, we have the following lemma.
\begin{lemma} \label{lem-operatorT}
The operator $I - \mathcal{T}_\eta$ is invertible 
with a bounded inverse in $X$.
\end{lemma}
Let $(E_\eta^h, H_\eta^h)$ be the solution to 
$$
(I-\mathcal{T}_\eta ) \begin{pmatrix}
E_\eta^{h} \\ H_\eta^{h}
\end{pmatrix}
=\begin{pmatrix}
E^{in} \\ H^{in}
\end{pmatrix}.
$$
Next, assume for simplicity that $\tilde{\epsilon}^N_{\mathrm{eff}}, \tilde{\mu}^N_{\mathrm{eff}},  \tilde{\tilde{\epsilon}}^N_{\mathrm{eff}},$ and
$\tilde{\tilde{\mu}}^N_{\mathrm{eff}}$ are replaced with their limits as  $N \rightarrow +\infty$ and 
consider the regularized form of (\ref{eq*}), that is,  
\begin{equation} \label{eq*eta}
\begin{pmatrix}
E^N_\eta(x)
\\
H^N_\eta(x)
\end{pmatrix}
= \begin{pmatrix}
E^{in}(x)
\\
H^{in}(x)
\end{pmatrix} + 
\frac{1}{N^3}\sum\limits_{j=1}^{N^3}  \omega G_\eta(x-z_j^N)
\begin{pmatrix}
\tilde{\epsilon}_{\mathrm{eff}} & i\omega\tilde{\tilde{\mu}}_{\mathrm{eff}} \\
-i\omega\tilde{\tilde{\epsilon}}_{\mathrm{eff}} & \tilde{\mu}_{\mathrm{eff}}
\end{pmatrix}
\begin{pmatrix}
E^{N}_\eta(z_j^N) \\ H^{N}_\eta(z_j^N)
\end{pmatrix}
\end{equation}
for $x\in\widetilde \Omega$. Here, $G_\eta$ is obtained from $G$ by replacing  $g^k$ with 
$g_{\eta}^k$,  and $(E^{N}_\eta(z_j^N), H^{N}_\eta(z_j^N))$ is obtained by solving the linear system
\begin{equation} \label{eq*eta2}
\begin{pmatrix}
E^N_\eta(z_i^N)
\\
H^N_\eta(z_i^N)
\end{pmatrix}
= \begin{pmatrix}
E^{in}(z_i^N)
\\
H^{in}(z_i^N)
\end{pmatrix} + 
\frac{1}{N^3}\sum\limits_{j=1, j\neq i}^{N^3}  \omega G_\eta(z_i^N -z_j^N)
\begin{pmatrix}
\tilde{\epsilon}_{\mathrm{eff}} & i\omega\tilde{\tilde{\mu}}_{\mathrm{eff}} \\
-i\omega\tilde{\tilde{\epsilon}}_{\mathrm{eff}} & \tilde{\mu}_{\mathrm{eff}}
\end{pmatrix}
\begin{pmatrix}
E^{N}_\eta(z_j^N) \\ H^{N}_\eta(z_j^N)
\end{pmatrix} , 
\end{equation}
for $i=1,\ldots, N^3$. 

To insure the uniform invertibility of (\ref{eq*eta2}) with respect to $N$ and $\eta$ (at least for $\omega$ small enough), we need some more assumptions regarding $\{z^N_j\}$ in addition to Assumption  \ref{dilute2}. We assume that 
$$
\frac{1}{N^6} \ds \sum_{i,j=1, i\neq j}^{N^3}  \| G(z_i^N - z_l^N) \|^2  \leq C N^{-6a}, 
$$
for some positive constant $C$. 

Define $\tilde{e}_\eta^N$  and  $\tilde{h}_\eta^N$    by $\tilde{e}_\eta^N:= E^h_\eta  - E_\eta^N$ and  
 $\tilde{h}_\eta^N:= H^h_\eta  - H_\eta^N$. Then, we have  
$$
\begin{array}{l}
\ds
\begin{pmatrix}
\tilde{e}^N_\eta(x)
\\
\tilde{h}^N_\eta(x)
\end{pmatrix}
-  
\frac{1}{N^3}\sum\limits_{j=1}^{N^3}  \omega G_\eta(x-z_j^N)
\begin{pmatrix}
\tilde{\epsilon}_{\mathrm{eff}} & i\omega\tilde{\tilde{\mu}}_{\mathrm{eff}} \\
-i\omega\tilde{\tilde{\epsilon}}_{\mathrm{eff}} & \tilde{\mu}_{\mathrm{eff}}
\end{pmatrix}
\begin{pmatrix}
\tilde{e}^{N}_\eta(z_j^N) \\ \tilde{h}^{N}_\eta(z_j^N)
\end{pmatrix} \\
\nm \ds = 
\omega \int_{{\Omega}}G_\eta(x-y)
\begin{pmatrix}
\tilde{\epsilon}_{\mathrm{eff}}& i\omega\tilde{\tilde{\mu}}_{\mathrm{eff}} \\
-i\omega\tilde{\tilde{\epsilon}}_{\mathrm{eff}} & \tilde{\mu}_{\mathrm{eff}}
\end{pmatrix}
\begin{pmatrix}
E^h_\eta(y)\\ H^h_\eta(y)
\end{pmatrix} \, dy - \frac{1}{N^3}\sum\limits_{j=1}^{N^3}  \omega G_\eta(x-z_j^N)
\begin{pmatrix}
\tilde{\epsilon}_{\mathrm{eff}} & i\omega\tilde{\tilde{\mu}}_{\mathrm{eff}} \\
-i\omega\tilde{\tilde{\epsilon}}_{\mathrm{eff}} & \tilde{\mu}_{\mathrm{eff}}
\end{pmatrix}
\begin{pmatrix}
E^{h}_\eta(z_j^N) \\ H^{h}_\eta(z_j^N)
\end{pmatrix}.
\end{array}
$$

Therefore, by the same arguments as those in \cite{acoustic, haisima}, 
one can prove that
$$
\|E_{\eta}^N - E^h_\eta\|_{\mathcal{C}^{0,\alpha}({\Omega})^3 } + \|H_{\eta}^N - H^h_\eta\|_{\mathcal{C}^{0,\alpha}({\Omega})^3 } \rightarrow 0 \quad \mbox{as } N \rightarrow +\infty,
$$
 uniformly in $\eta$. 
Then,  since on one hand, $E_{\eta}^N \rightarrow E^{N}$ and $H_{\eta}^N \rightarrow H^{N}$ in ${\Omega} \setminus \cup (z_j)_{j=1}^{N^3}$ as $\eta\rightarrow0$ and on the other hand, 
$E_{\eta}^h \rightarrow E^{h}$ and $H_{\eta}^h \rightarrow H^{h}$ in ${\Omega}$ as $\eta\rightarrow0$, 
we obtain the desired justification of the homogenization procedure.

\end{document}